\documentclass[12pt,a4paper]{article}

\usepackage[reqno]{amsmath}
\usepackage{amssymb,amsthm,upref,amscd}
\usepackage[T1]{fontenc}
\usepackage{times}
\usepackage{color}

\newtheorem{theorem}{Theorem}[section]
\newtheorem{corollary}[theorem]{Corollary}
\newtheorem{definition}[theorem]{Definition}

\newtheorem{lemma}[theorem]{Lemma}

\newtheorem{remark}[theorem]{Remark}

\theoremstyle{definition} \theoremstyle{remark}
\numberwithin{equation}{section}

\begin{document}

\title{\textbf{Uniqueness of normalized ground states for NLS models  \ \\
}}
\author{Hichem Hajaiej$^{\mathrm{a,}}$, Linjie Song$^{\mathrm{b,}}$\thanks{%
Linjie Song is supported by "Shuimu Tsinghua Scholar Program" and by "National Funded Postdoctoral Researcher Program" (GZB20230368), and funded by "China Postdoctoral Science Foundation" (2024T170452). Email: songlinjie18@mails.ucas.edu.cn.}\  \\
\\
{\small $^{\mathrm{a}}$ Department of Mathematics, California State University at Los Angeles, Los Angeles, CA 90032, USA}\\
{\small $^{\mathrm{b}}$Department of Mathematical Sciences, Tsinghua University, Beijing 100084, China}
}
\date{}
\maketitle

\begin{abstract}
	
We present two methods to prove 
the uniqueness of normalized ground states. We will first discuss the key ideas and ingredients of each method. Then, we will apply them to various classes of PDEs. Our approach is applicable to other operators, domains and nonlinearities provided that some hypotheses are satisfied. 

\bigskip



\end{abstract}

\medskip

\section{Introduction} \label{int}

In \cite{HS}, we provided an abstract framework to prove the uniqueness of the ground states on Nehari manifolds. This has enabled us to derive results on the existence and non-existence of the \emph{normalized solutions} for a large class of PDEs (see applications in \cite{LS,Song,Song2,Song3,Song4}). Normalized solutions, that are solutions with prescribed $L^2$ mass $c > 0$, have gained great interest in the last few years due to their numerous applications in physics. In particular, they play an important role in studying the orbital stability of standing waves. The main goal of this work is to discuss the uniqueness of this important class of solutions by providing two general approaches. We will also connect this property to the strict monotonicity of the global branch of solutions of the original PDE.  More precisely, we consider a general equation given by:
\begin{equation} \label{eqabs}
	D_u\Phi_\lambda(u) = 0, \quad u \in W,
\end{equation}
where $W$ is a Banach space with the norm $\|\cdot\|$, $\Phi_{\lambda} \in C^{2}(W, \mathbb{R})$ is a Fr\'{e}chet-differentiable functional with derivative $D_{u}\Phi_{\lambda}: W \rightarrow W^{\ast}$, the space $W^{\ast}$ is the dual space of $W$, $\lambda \in \mathbb{R}$. In the first part of our abstract framework developed in \cite{HS}, we established general conditions under which, we have the existence of a $C^1$ global branch
$$
u: (-\infty,\lambda_1) \rightarrow W,
$$
such that $u(\lambda)$ solves \eqref{eqabs} and has Morse index $1$, where $\lambda_1$ is some real number. Furthermore, we showed that any solution with Morse index $1$ is unique and lies on this branch. (The general setting and results on the existence and uniqueness of this global branch are shown in Appendix \ref{exi and uni}.)
From now on, we will assume that $\Phi_\lambda$ has the form 
\begin{align}
	\Phi_\lambda(u) = E(u) - \lambda Q(u) \tag*{(1.1)$_E$},
\end{align}
where 
\begin{align}
2Q(u) = D_uQ(u)(u). \tag*{(1.1)$_Q$}
\end{align}
If $E$ is bounded from below on $S_c := \{u \in W: Q(u) = c\}$, the minimizers under the $L^2$ constraint are the functions achieving  $m(c) := \inf_{u \in S_c}E(u)$. In this case, $\lambda$ is a Lagrange multiplier. Such solutions are called \emph{normalized ground states} on $S_c$ in this paper. They are known to be the best candidates to enjoy stability. The existence of normalized ground states has been widely studied. However the uniqueness has only been addressed for simple nonlinearities -typically pure power nonlinearities-. The main difficulty of the establishment of this important property is that several Lagrange multipliers can give the same mass $Q$ and energy $E$. This means that the uniqueness of the solutions to the equation at a fixed Lagrange multiplier does not imply the uniqueness of the normalized ground solutions. On the contrary, the existence (respectively the uniqueness) of the normalized ground solutions implies the existence (respectively the uniqueness) of the minimizers on Nehari manifold at a fixed Lagrange multiplier under certain conditions, see the proof in a very general setting in Appendix \ref{Neh}. Hence, the uniqueness of the normalized ground solutions is much more challenging, which would explain the silence of the literature. 
Many authors conjectured results about the uniqueness of normalized ground states (see remarks, page 675, \cite{SS}, Abstract and Introduction of \cite{LN}). Others used numerical simulations to convince the community of the validity of this uniqueness. Let us point out that \cite{LN} treated the delicate cubic-quintic case and brought significant contribution to the topic.

The uniqueness of the normalized ground states is heavily connected to the strong orbital stability of standing waves that was introduced by Grillakis, Shatah, and Strauss in their breakthrough paper \cite{GSS1}. Their main result \cite[Theorem 3.2]{GSS1} was a source of inspiration for many mathematicians and physicists especially that their assumptions seem hard to be checked. The aim of this paper is not only to offer an alternative to it but to also present two general approaches to prove the uniqueness of the normalized ground states. More precisely, we present two general, flexible and adaptable methods to prove the strict monotonicity of the branch of solutions, and to show the uniqueness of the corresponding normalized ground states.

Let $u_\lambda$ be a solution of $\eqref{eqabs}_E$, we now introduce two main paths to achieve this goal:

\vskip0.1in
\underline{Method 1:} There are three main steps in this approach:

\vskip0.1in
\noindent\underline{Step 1:} Prove that $v_\lambda=\partial_\lambda u_\lambda$ changes sign at most once in $r>0$.\\
\underline{Step 2:} Show that it is impossible to find a $\lambda<\lambda_1$ such that $\partial_\lambda |u_\lambda|_2=0$, $|\cdot|_2$ is the $L^2-$norm.\\
\underline{Step 3:} By step 2, $\partial_\lambda|u_\lambda|_2>0,$ or $\partial_\lambda|u_\lambda|_2<0,$ for all $\lambda<\lambda_1.$ Using the asymptotic behavior of $|u_\lambda|_2$ when $\lambda\to -\infty$ or $\lambda\to \lambda_1$ enables us to conclude.

\vskip0.1in
In our applications, step 1 plays an important role in showing step 2. If $\partial_\lambda |u_\lambda|_2=0$ for some $\lambda$, $v_\lambda$ is sign-changing. By step 1, $v_\lambda$ exactly changes sign once. Then we can derive a contradiction. In order to prove step 1, we will show the following properties
$$
(\text{H}1)\quad u_\lambda^{'}(r)<0\ \text{for}\ r>0;
$$
$$
(\text{H}2)\quad v_\lambda(0)=\frac{d }{d \lambda}u_\lambda(0)<0.
$$

Taking into account a combination of the local and global properties of the solutions, we were able to develop a general approach by providing all the details in two concrete examples. The first example, Eq \eqref{eq1.1}, has been playing a crucial role in nonlinear optics. Note that $(\text{H}1)$ and $(\text{H}2)$ imply that $v_\lambda=\frac{d }{d \lambda}u_\lambda$ changes sign at most once in $r$ (see Lemma \ref{lemB.1}). The latter is critical to reach our goal. $(\text{H}2)$ can be replaced/weakened when dealing with PDEs on bounded domains. For the class of PDEs given by \eqref{eq1.2}, we will assume $(\text{H}2^{'}):$ $v_\lambda(r)<0$ in the neighborhood of a certain $r$ when $\frac{d}{d \lambda}|v_\lambda|_2=0$. $(\text{H}1)$ and $(\text{H}2^{'})$ also imply that $v_\lambda$ changes sign exactly ones (see Proof of Lemma \ref{lemD.2}). This approach seems more appropriate when one deals with situations on bounded domains. Let us point out that we believe that this method is only valid for local operators. However, it has the considerable advantage of not requiring the boundedness of the energy from below. Therefore, it can address all the regimes (subcritical, critical, and supercritical).

 We now introduce another new and self-contained approach that will enable us to prove the uniqueness of normalized ground states.

 \vskip0.1in
\underline{Method 2:} There are two major steps to use this approach:

\vskip0.1in
\noindent\underline{Step 1:} Show that $m(c)$ is differentiable at $c$ if and only if
 $$
\Gamma(c) = \big\{\lambda: \exists u_\lambda \ \text{satisfying} \ \eqref{eqabs} \ \text{such that} \ Q(u_\lambda) = c, E(u_\lambda) = m(c)\big\}.
$$
has exactly one element.

\noindent\underline{Step 2:} Prove that $\Gamma(c)$ has exactly one element if and only if the normalized ground state on $S_c$ is unique.

  \vskip0.1in
Step 1 and Step 2 show that:

{\bf  The normalized ground state on $S_c$ is unique $\Leftrightarrow$ $m$ is differentiable at $c.$}


We believe that Method 2 is complementary to the first method. The smoothness of $m(c)$ has been touched upon in \cite{Stefanov}. Method 2 is the first result in the literature binding, in an equivalent manner, the uniqueness of the normalized ground state on $S_c$ and the differentiability of the function $m$ at $c$. 

Now, to prove step 2, it is sufficient to show that all the normalized ground states have Morse index $1,$ and that the solutions to \eqref{eqabs} with Morse index $1$ is unique at a fixed $\lambda$. On the other hand, the differentiability of $m$ almost everywhere is guaranteed by its monotonicity.

Contrary to Method 1, this method requires the boundedness from below of the energy functional $E$ on the sphere.


For the readers' convenience, we study concrete equations to show how our general framework works. Our approach is applicable to many other equations. Let us focus for now on a class of PDEs that has attracted the attention of several colleagues when $h\equiv1$. Namely:
\begin{equation} \label{eq1.4}
	\begin{cases}
		(-\Delta)^s u  = \lambda u + h(|x|)|u|^{p-2}u
		& \text{in } \mathbb{R}^N, \\[1.5\jot]
		\displaystyle
		u(x) \rightarrow 0
		&\text{as } |x| \rightarrow +\infty,
	  \end{cases}
\end{equation}
where $0 < s \leq 1$. When $s = 1$, $(-\Delta)^s$ is the usual Laplacian operator while when $s \in (0,1)$, $(-\Delta)^s$ is the fractional Laplacian, see Section 5 in \cite{HS} for more details about its definition. We assume that $h$ satisfies the following assumptions:

\begin{itemize}
	\item[$(h)$] $h(r) \in C^{1}([0,+\infty)) \cap L^{\infty}([0,+\infty))$, $h(r) > 0$ in $[0,+\infty)$, $h(r)$ and $\frac{rh'(r)}{h(r)}$ are non-increasing in $(0,+\infty)$, $\theta = \lim_{r \rightarrow +\infty}\frac{rh'(r)}{h(r)} > -2s$, $\sup_{r > 0}r^{-\theta}h(r) < \infty$.
\end{itemize}

In this case,
\begin{align*}
	\Phi_\lambda(u) & = \frac{1}{2}\int_{\mathbb{R}^N}\left(|(-\Delta)^{\frac{s}{2}} u|^2 - \lambda |u|^2\right)dx - \frac{1}{p}\int_{\mathbb{R}^N}h(|x|)|u|^{p}dx, \\
E(u) & = \frac{1}{2}\int_{\mathbb{R}^N}|(-\Delta)^{\frac{s}{2}} u|^2dx - \frac{1}{p}\int_{\mathbb{R}^N}h(|x|)|u|^{p}dx, \\
Q(u) & = \frac{1}{2}\int_{\mathbb{R}^N}|u|^2dx, \\
m(c) & = \inf_{u \in S_c}E(u) \quad \text{where } S_c = \big\{u \in H^s(\mathbb{R}^N): Q(u) = c \big\}.
\end{align*}
Let
$$
\lambda_{1,c} := \inf\Lambda(c), \quad \lambda_{2,c} := \sup\Lambda(c),
$$
where
$$
\Lambda(c) := \big\{\lambda: \exists u_\lambda \text{ satisfying } \eqref{eq1.4} \text{ such that } Q(u_\lambda) = c, E(u_\lambda) = m(c)\big\}.
$$
In this context, Method 2 will enable us to show that:

\begin{theorem}
	\label{thm2.10}	Assume that $(h)$ holds and $2 < p < 2 + (2\theta + 4s)/N$. Then
	\begin{itemize}
		\item[$(i)$] $m(c)$ is differentiable at almost every $c > 0$.
		\item[$(ii)$] We have $m'_{+}(c) = \lambda_{1,c}$ and $m'_{-}(c) = \lambda_{2,c}$. In particular, $m(c)$ is differentiable at $c$ if and only if $\Lambda(c)$ has only one element.
		\item[$(iii)$] For almost every $c > 0$, \eqref{eq1.4} has a unique normalized ground state (after multiplying $-1$ if necessary) on $S_c$, which is positive.
	\end{itemize}
\end{theorem}

\begin{remark}
	We work in a setting where the branch of solutions is unique. In this case, there is equivalence between the strict decreasiness of the $L^2$ norm of the solutions on the global branch and the uniqueness of the normalized ground states for all $c > 0$.
\end{remark}

\begin{remark}
	Let $u_c$ be a positive normalized ground state on $S_c$ for \eqref{eq1.4} with Lagrange multiplier $\lambda$, and let $\Psi(t,x) = e^{-i\lambda t}u_c(x)$. Then $\Psi(t,x)$ is a standing wave solution of the following Cauchy problem:
\begin{equation} \label{cauchypro}
	\begin{cases}
		i\partial_t\Psi - (- \Delta)^s \Psi =  -h(|x|)|\Psi|^{p-2}\Psi \quad \text{in } \mathbb{R}_+ \times \mathbb{R}^N, \\[1.5\jot]
		\displaystyle
		\Psi(0,x) = \Psi_0(x),
	  \end{cases}
\end{equation}
with $\Psi_0(x) = u_c(x)$. Note that $h \in L^\infty$ and $2 < p < 2 + (2\theta + 4s)/N < 2 + 4s/N$. It is well known that the Cauchy problem \eqref{cauchypro} is globally well-posed in $H^s$, see \cite[Corollary 6.1.2]{Caz} for $s = 1$ and see \cite[Theorem 2.6]{GH} for $s \in (0,1)$. Let 
\begin{align*}
	\mathcal{M}_c := \big\{e^{i\theta}u: \theta \in \mathbb{R}, u \text{ is a normalized ground state of } \eqref{eq1.4} \text{ on } S_c\big\}.
\end{align*}
Following the arguments in \cite{CL} (see \cite{Al} for an application to nonlocal nonlinear wave equations), one can verify that $\mathcal{M}_c$ is stable by the flow of \eqref{cauchypro} (see the definition and more details in \cite{CL}, also called weak stability in some articles). When the positive normalized ground state is unique, the weak stability is equivalent to the strong (also called true) orbital stability (see the definition and more discussions in \cite[p.370]{Stuart}).
\end{remark}


\begin{remark}
	\begin{itemize}
		\item[$(i)$] We construct a counter-example where the normalized ground state exists but is not unique (after multiplying $-1$ if necessary) in Appendix A.
		\item[$(ii)$] If the normalized ground state is unique for $c > 0$, we give the explicit expression of $m(c)$ in terms of $u_c,$ the unique normalized ground state solution in Appendix D. If $m$ is differentiable at $c$, Theorem \ref{thm2.10} states that there is a unique normalized ground solution. Again Appendix D provides the expression of $m$. This expression was only known for very particular nonlinearities ($s=1,$ $h(|x|)=|x|^\alpha,$ $\alpha>0$, or pure power nonlinearity \cite{MM}, $0<s<1,$ $h(|x|)\equiv1,$ \cite{Luo-Zhang}).
	\end{itemize}
\end{remark}

\begin{remark}
	Without discussing the details, below are other examples that our method applies to:
	\begin{itemize}
		\item[$(i)$] The equation involving the mixed fractional Laplacian, $(-\Delta)^{s_1} u + (-\Delta)^{s_2} u = \lambda u + |u|^{p-2}u$ where $s_1 < s_2 \leq 1 \leq 2s_1$, $2 < p < 2+4s_1/N$ and $s_2 - s_1$ is small enough (independently of $\lambda$).
		\item[$(ii)$] The cubic-quintic NLS, $-\Delta u = \lambda u + |u|^{2}u - |u|^{4}u$ in $\mathbb{R}^3$. We remark that this equation has been addressed in \cite{LN}. 
	\end{itemize}
\end{remark}

\begin{remark}
	We thank the anonymous referees for pointing out reference \cite{DST1}, where the authors studied the uniqueness and nonuniqueness of normalized ground states on metric graphs. The proof of Theorem \ref{thm2.10} is based mainly in Theorem \ref{thm5.1}, Lemma \ref{lem5.4} and similar ideas to the work of \cite{DST1}.
\end{remark}

\begin{remark}
	Including all the values of $c$ in Theorem \ref{thm2.10} $(iii)$ seems very challenging. For small and large $c,$ this is the subject of an ongoing work \cite{HPS}.
\end{remark}


We now provide two examples, in which we use Method 1 to study the strict monotonicity of the $L^2$ norm of $u(\lambda)$ with respect to $\lambda$  on the global branch for all $\lambda \in (-\infty,\lambda_1)$:
\begin{equation} \label{eq1.1}
	\begin{cases}
		-\Delta u + V(|x|)u = \lambda u + f(|x|,u)
		& \text{in } \mathbb{R}^N, \\[1.5\jot]
		\displaystyle
		u(x) \rightarrow 0
		&\text{as } |x| \rightarrow +\infty,
	  \end{cases}
\end{equation}
\begin{equation} \label{eq1.2}
	\begin{cases}
		- \Delta u = \lambda u + |x|^{-k}|u|^{p-2}u
		& \text{in } B_1, \\[1.5\jot]
		\displaystyle
		u(x) = 0
		&\text{on } \partial B_1,
	  \end{cases}
\end{equation}
$B_1 \subset \mathbb{R}^N$ is the unit ball with $N \geq 3$, and $0 < k < 2$. This implies the uniqueness of the normalized ground states. Moreover, this strict monotonicity has a close connection with the orbital stability of these solutions, when viewed as standing waves. Despite the importance of this question, there are very few results in the literature addressing this aspect. To the best of our knowledge, the only paper is \cite{MST}, where the authors studied this issue for equation \eqref{eq1.1} in the $L^2$ subcritical regime when $N = 1$.

For \eqref{eq1.1} in the $L^2$ subcritical case, let $u_\lambda = u(\lambda) \in H^1(\mathbb{R}^N)$, and $\lambda_1 =
\inf \sigma(-\Delta + V)$ be the first eigenvalue. Assuming that $u_\lambda(x) = u_\lambda(|x|)$, $v_\lambda(0) < 0$ and under some other suitable assumptions (see Subsection \ref{nls} for more details), we will show that $v_\lambda=\partial_\lambda u_\lambda$ changes sign at most once in $r > 0$. This is critical to show that it is impossible that
$$
\partial_\lambda\int_{\mathbb{R}^N}|u_\lambda|^2dx = 2\int_{\mathbb{R}^N}u_\lambda v_\lambda dx = 0, \quad \forall \lambda < \lambda_1.
$$
Then classical bifurcation arguments (c.f. \cite[Theorem 1.7]{CR}) yield that $\lim_{\lambda \rightarrow \lambda_1}\int_{\mathbb{R}^N}|u_\lambda|^2dx = 0$, implying that
$$
\partial_\lambda\int_{\mathbb{R}^N}|u_\lambda|^2dx < 0, \quad \forall \lambda < \lambda_1.
$$
We will further verify the hypothesis $v_\lambda(0) < 0$ when $N = 1$. Moreover, we will consider the case when $V \equiv 0$ using similar arguments. Both the $L^2$ subcritical case and the $L^2$ supercritical case will be addressed for \eqref{eq1.1}. When $V \equiv 0$, $\lambda_1 = 0$ in \eqref{eq1.1}, we will apply the method developed in \cite{HS} to determine the asymptotic behaviors of $\int_{\mathbb{R}^N}|u_\lambda|^2dx$ when $\lambda \rightarrow 0^-$ or $\lambda \rightarrow -\infty$ to complete the proof of $\partial_\lambda\int_{\mathbb{R}^N}|u_\lambda|^2dx < 0$ or $\partial_\lambda\int_{\mathbb{R}^N}|u_\lambda|^2dx > 0$ for all $\lambda < 0$.\\

For \eqref{eq1.2}, let $u_\lambda = u(\lambda) \in H_{0,rad}^1(B_1)$, and $\lambda_1 = \lambda_1(B_1)$ be the first eigenvalue of $-\Delta$ on $B_1$ with Dirichlet boundary condition. We no longer assume that $v_\lambda(0) < 0$ where $v_\lambda = \partial_\lambda u_\lambda$. Instead, we will prove that $v_\lambda(r) < 0$ near $r = 1$ if $\int_{B_1}u_\lambda v_\lambda dx = 0$. Then under suitable assumptions (see Subsection \ref{nlsb} for more details), similar arguments to \eqref{eq1.1}, will enable us to show that:
$$
\partial_\lambda\int_{\mathbb{R}^N}|u_\lambda|^2dx < 0, \quad \forall \lambda < \lambda_1(B_1),
$$
and that the normalized ground state is unique.

In Section \ref{pr}, we prove Theorem \ref{thm2.10} by showing all the steps described previously. In Section \ref{result}, we state the results for \eqref{eq1.1}, \eqref{eq1.2}, and give their proofs in Sections \ref{pr1}, \ref{pr2} respectively. As mentioned before, many other operators, nonlinearities, and domains can be included in our approach.

\section{Proof of Theorem \ref{thm2.10}} \label{pr}

\begin{theorem}[Theorem 5.16 and Lemma 5.13 in \cite{HS}]
	\label{thm5.1} Assume that $(h)$ holds and that $2 < p < 2 + (2\theta + 4s)/N$. Then for any $\lambda < 0$, \eqref{eqabs} has a unique positive solution $u_\lambda$ with Morse index $1$. Furthermore, $u_\lambda$ is radial and non-degenerate.
\end{theorem}

\begin{lemma}
	\label{lem 5.2} Assume that $(h)$ holds and $2 < p < 2 + (2\theta + 4s)/N$. Then $-\infty < m(c) < 0$.
\end{lemma}

\begin{proof} The fractional Gagliardo-Nirenberg inequality (see \cite{HH1}) shows that
\begin{equation}
\int_{\mathbb{R}^N}|u|^{p}dx \leq C(s,N,p)\left( \int_{\mathbb{R}^N}|(-\Delta)^{\frac{s}{2}} u|^2dx\right)^{\frac{N(p-2)}{4s}}\left( \int_{\mathbb{R}^N}|u|^2dx\right)^{\frac{p}{2}-\frac{N(p-2)}{4s}}.
\end{equation}
Hence, for any $u \in S_c$, we have
\begin{equation} \label{eq5.2}
\int_{\mathbb{R}^N}h(|x|)|u|^{p}dx \leq C(s,N,p)\|h\|_{L^{\infty}}\left( 2c\right)^{\frac{p}{2}-\frac{N(p-2)}{4s}}\left( \int_{\mathbb{R}^N}|(-\Delta)^{\frac{s}{2}} u|^2dx\right)^{\frac{N(p-2)}{4s}}.
\end{equation}
Note that
$$
0 < \frac{N(p-2)}{4s} < 1.
$$
We get that
$$
m(c) = \inf_{u \in S_c}E(u) > -\infty.
$$

Now we aim to prove that $m(c) < 0$. To do this, we first show that $h(tr) \geq \epsilon \kappa(t,r)$ for all $t > 0, r > 0$ with some $\epsilon > 0$ where $\kappa(t,r) = t^{\theta}h(r)$ for $tr > 1$ and $\kappa(t,r) = h(r)$ for $tr \leq 1$. Indeed, on the one hand, let us consider
$$
g(t) = h(tr) - \epsilon_1t^{\theta}h(r), \quad t > 1/r.
$$
Note that
$$
g'(t) = h'(tr)r-\epsilon_1\theta t^{\theta-1}h(r) \geq \frac{\theta}{t}g(t).
$$
Since $\sup_{r>0}r^{-\theta}h(r) < \infty$, we can take $\epsilon_1 > 0$ such that $g(\frac1r) > 0$. Then by Gronwall inequality, one gets that $g(r) > 0$ for all $t > 1/r$. On the other hand, for $t \leq 1/r$,
$$
h(tr) \geq h(1) \geq \frac{h(1)}{\|h\|_{L^\infty}}h(r).
$$
Hence, take $\epsilon = \min\{\epsilon_1,\frac{h(1)}{\|h\|_{L^\infty}}\}$ and we obtain that $h(tr) \geq \epsilon \kappa(t,r)$.

For any $u \in S_c$, let
$$
u_\tau = \tau^{\frac{N}{2}}u(\tau x) \in S_c.
$$
Then for $\tau < 1$,
\begin{eqnarray}
&& E(u_\tau) \nonumber \\
&=& \frac{1}{2}\tau^{2s}\int_{\mathbb{R}^N}|(-\Delta)^{\frac{s}{2}} u|^2dx - \frac{1}{p}\int_{\mathbb{R}^N}h(|x|)|u_\tau|^{p}dx \nonumber \\
&\leq& \frac{1}{2}\tau^{2s}\int_{\mathbb{R}^N}|(-\Delta)^{\frac{s}{2}} u|^2dx - \frac{\epsilon}{p}\left(\tau^{\frac{p-2}{2}N-\theta}\int_{|x|>\tau}h(|x|)|u|^{p}dx + \tau^{\frac{p-2}{2}N}\int_{|x|\leq\tau}h(|x|)|u|^{p}dx\right) \nonumber \\
&\leq& \frac{1}{2}\tau^{2s}\int_{\mathbb{R}^N}|(-\Delta)^{\frac{s}{2}} u|^2dx - \frac{\epsilon}{p}\tau^{\frac{p-2}{2}N-\theta}\int_{\mathbb{R}^N}h(|x|)|u|^{p}dx.
\end{eqnarray}
Since $2s > (p-2)N/2-\theta$, $E(u_\tau) < 0$ if $\tau$ is small enough, implying that $m(c) < 0$.
\end{proof}

\begin{lemma}
	\label{lem 5.3}  Assume that $(h)$ holds and $2 < p < 2 + (2\theta + 4s)/N$. Let $c_1, c_2 > 0$ be such that $c = c_1 + c_2$. Then
	$$
	m(c) < m(c_1) + m(c_2).
	$$
\end{lemma}

\begin{proof} 
	Let $\{u_n\} \subset S_c$ be a minimizing sequence for $m(c)$. For any $k > 1$, we have
\begin{equation} \label{eqk2c}
	m(k^2c) \leq E(ku_n) = \frac{k^{2}}{2}\int_{\mathbb{R}^N}|(-\Delta)^{\frac{s}{2}} u_n|^2dx - \frac{k^{p}}{p}\int_{\mathbb{R}^N}h(|x|)|u_n|^{p}dx < k^{2}E(u_n) \to k^2m(c).
\end{equation}
 Thus $m(k^2c) \leq k^2m(c)$. We claim that $m(k^2c) < k^2m(c)$. Arguing by contradiction, we assume that $m(k^2c) = k^2m(c)$. Then \eqref{eqk2c} indicates
 \begin{align}
 	& \lim_{n \to \infty}\left( \frac{k^{2}}{2}\int_{\mathbb{R}^N}|(-\Delta)^{\frac{s}{2}} u_n|^2dx - \frac{k^{p}}{p}\int_{\mathbb{R}^N}h(|x|)|u_n|^{p}dx\right) \nonumber \\
 	= \ & \lim_{n \to \infty}\left( \frac{k^{2}}{2}\int_{\mathbb{R}^N}|(-\Delta)^{\frac{s}{2}} u_n|^2dx - \frac{k^{2}}{p}\int_{\mathbb{R}^N}h(|x|)|u_n|^{p}dx\right).
 \end{align}
Consequently, $\int_{\mathbb{R}^N}h(|x|)|u_n|^{p}dx \rightarrow 0$ as $n \to \infty$, showing that
\begin{equation}
\lim_{n \rightarrow \infty}E(u_n) = \frac{1}{2}\lim_{n \rightarrow \infty}\int_{\mathbb{R}^N}|(-\Delta)^{\frac{s}{2}} u_n|^2dx \geq 0,
\end{equation}
contradicting to the fact that $m(c) < 0$. Hence, the claim is proved.

Now we prove that $m(c) < m(c_1) + m(c_2)$. WLOG, we may assume that $c_1 \leq c_2$. Then
\begin{equation}
m(c) = m(\frac{c}{c_2}c_2) < \frac{c}{c_2}m(c_2) = m(c_2) + \frac{c_1}{c_2}m(c_2) = m(c_2) + \frac{c_1}{c_2}m(\frac{c_2}{c_1}c_1) \leq m(c_2) + m(c_1).
\end{equation}
The proof is now complete.
\end{proof}

\begin{corollary}
	\label{cor5.4} Assume that $(h)$ holds and $2 < p < 2 + (2\theta + 4s)/N$. Then $m(c)$ is strictly decreasing with respect to $c > 0$.
\end{corollary}

\begin{proof}
	Let $0 < c_1 < c$, we get $m(c) < m(c_1)$ immediately. In fact, by Lemma \ref{lem 5.2} and Lemma \ref{lem 5.3}, we get
	\begin{align*}
		m(c) < m(c_1) + m(c-c_1) < m(c_1).
	\end{align*}
	The proof is complete.
\end{proof}

Like the proof of \cite[Lemma 4.2]{YQZ} or \cite[Theorem 3.3]{EHLS}, we have the following lemma and omit the proof here.

\begin{lemma}
 Assume that $(h)$ holds and $2 < p < 2 + (2\theta + 4s)/N$. Then $m(c)$ is continuous with respect to $c > 0$.
\end{lemma}

\begin{lemma}
	\label{lem5.2} Assume that $(h)$ holds and $2 < p < 2 + (2\theta + 4s)/N$. Let $\{u_n\}$ be a sequence in $H_{rad}^s$ such that
	$$
	E(u_n) \rightarrow m(c) \text{ and } Q(u_n) \rightarrow c.
	$$
	If $N \geq 2$ or if $u_n(x)$ is a nonincreasing function of $|x|$ for every $n$, then there is a subsequence, also denoted by $\{u_n\}$, and $u \in H^s$ such that $u_n \to u$ strongly in $H^s$ as $n \to \infty$. As a consequence, for each $c > 0$, there exists a $u_c \in S_c$ such that $E(u_c) = m(c)$. Furthermore, any normalized ground state on $S_c$ is positive (after multiplying $-1$ if necessary) and radially symmetric.
\end{lemma}

\begin{proof} 
	Let $\{u_n\}$ be a sequence in $H_{rad}^s$ such that $E(u_n) \rightarrow m(c)$ and $Q(u_n) \rightarrow c$. By \eqref{eq5.2}, we have
	\begin{align}
		E(u_n) \geq & \ \frac 12\int_{\mathbb{R}^N}|(-\Delta)^{\frac{s}{2}} u_n|^2dx \\
		& - \frac{C(s,N,p)}p\left( \int_{\mathbb{R}^N}|(-\Delta)^{\frac{s}{2}} u_n|^2dx\right)^{\frac{N(p-2)}{4s}}\left( \int_{\mathbb{R}^N}|u_n|^2dx\right)^{\frac{p}{2}-\frac{N(p-2)}{4s}}.
	\end{align}
	Since $N(p-2)/4s < 1$ , we get $\int_{\mathbb{R}^N}|(-\Delta)^{\frac{s}{2}} u_n|^2dx$ is bounded. Therefore, $\{u_n\}$ is bounded in $H^s$. Up to a subsequence, we may assume that there exists $u \in H_{rad}^s$ such that $u_n \rightharpoonup u$ weakly in $H^s$ and $u_n \rightharpoonup u$ weakly in $L^2$. By \cite[Proposition 1.7.1]{Caz}, we also assume that $u_n \rightarrow u$ strongly in $L^p$. Note that $h(r)$ is bounded. Then, by the weak convergence of $u_n$ to $u$ in $H^s$ and strong convergence in $L^p$, we get
$$
E(u) \leq \lim_{n \rightarrow \infty}E(u_n) = m(c),
$$
Moreover, by the weak convergence of $u_n$ to $u$ in $L^2$, we have
$$
Q(u) \leq \lim_{n \rightarrow \infty}Q(u_n) = c.
$$
Next, by the definition of $m(c)$ and Corollary \ref{cor5.4}, we obtain
\begin{align*}
	m(Q(u)) \leq E(u) \leq m(c) \leq m(Q(u)),
\end{align*}
which implies that $E(u) = m(c)$ and $Q(u) = c$. The weak convergence of $u_n$ in $L^2$ and the fact that $Q(u_n) \rightarrow Q(u)$ show that $u_n$ strongly converges to $u$ in $L^2$. Similarly, we can show the strong convergence of $u_n$ in $H^s$.


Next, we show the existence of a normalized ground state. Let $\{u_n\} \subset S_c$ be a minimizing sequence for $m(c)$. Using the Schwarz rearrangement of $u_n$ (see \cite{HH2} or \cite{Kaw}), we may assume that $u_n$ is nonnegative, radially symmetric and nonincreasing of $|x|$. Obviously, $Q(u_n) = c$ implies that $Q(u_n) \rightarrow c$ holds. Hence, we obtain the existence of a minimizer $u_c$ immediately. The Schwarz rearrangement of $u_c$ shows that it is non-negative (after multiplying $-1$ if necessary) and radially symmetric. Positivity is a direct consequence of the maximum principle. The proof is complete.
\end{proof}

\begin{lemma}
	\label{lem5.4} Assume that $(h)$ holds and $2 < p < 2 + (2\theta + 4s)/N$. Then for each $c > 0$, we have
	\begin{itemize}
		\item[$(i)$] $-\infty < \inf \Lambda(c) \leq \sup \Lambda(c) < 0$;
		\item[$(ii)$] $\lambda_{1,c}, \lambda_{2,c} \in \Lambda(c)$.
	\end{itemize} 
\end{lemma}

\begin{proof} First, we show that $-\infty < \inf \Lambda(c) \leq \sup \Lambda(c) \leq 0$. If $u \in S_{c}$ solves \eqref{eq1.4}, then the following integral identity (c.f. \cite[Lemma 5.6]{HS}) holds true:
\begin{equation}
(N-2s)\int_{\mathbb{R}^N}|(-\Delta)^{\frac{s}{2}} u|^2dx = N\lambda \int_{\mathbb{R}^N}|u|^2dx + \frac{2}{p}\int_{\mathbb{R}^N}(Nh(|x|)+|x|h'(|x|))|u|^{p}dx.
\end{equation}
Note that
\begin{equation}
\int_{\mathbb{R}^N}|(-\Delta)^{\frac{s}{2}} u|^2dx = \lambda \int_{\mathbb{R}^N}|u|^2dx + \int_{\mathbb{R}^N}h(|x|)|u|^{p}dx.
\end{equation}
Therefore, we can obtain
\begin{equation} \label{eq5.5}
-2s\lambda \int_{\mathbb{R}^N}|u|^2dx = \int_{\mathbb{R}^N}\left( \left(\frac{2N}{p}-(N-2s)\right) h(|x|)+\frac{2}{p}|x|h'(|x|)\right) |u|^{p}dx,
\end{equation}
\begin{equation} \label{eq5.6}
2s\int_{\mathbb{R}^N}|(-\Delta)^{\frac{s}{2}} u|^2dx = \int_{\mathbb{R}^N}\left( \frac{p-2}{p}N h(|x|)-\frac{2}{p}|x|h'(|x|)\right) |u|^{p}dx.
\end{equation}

On the one hand, if $\lambda \rightarrow -\infty$, \eqref{eq5.5} shows that
$$
\int_{\mathbb{R}^N}h(|x|)|u|^{p}dx \rightarrow +\infty.
$$
Then from \eqref{eq5.6} and $(h)$ we deduce that
\begin{eqnarray}
E(u) &=& \int_{\mathbb{R}^N}\left( \left( \frac{p-2}{4sp}N - \frac{1}{p}\right) h(|x|)-\frac{1}{2sp}|x|h'(|x|)\right) |u|^{p}dx \nonumber \\
&\leq& \left( \frac{p-2}{4sp}N - \frac{1}{p}-\frac{\theta}{2sp}\right)\int_{\mathbb{R}^N}h(|x|)|u|^{p}dx \nonumber \\
&\rightarrow& -\infty,
\end{eqnarray}
since $p < 2 + (2\theta + 4s)/N$ implies that 
$$
\frac{p-2}{4sp}N - \frac{1}{p}-\frac{\theta}{2sp} < 0.
$$
However, by Lemma \ref{lem 5.2},
$$
m(c) = \inf_{u \in S_c}E(u) > -\infty,
$$
which is a contradiction. Thus $\inf \Lambda(c)  > -\infty$.

On the other hand, by \eqref{eq5.5}, we obtain that
\begin{equation} \label{eq5.8}
-2s\lambda \int_{\mathbb{R}^N}|u|^2dx \geq \int_{\mathbb{R}^N} \left(\frac{2(N+\theta)}{p}-(N-2s)\right) h(|x|) |u|^{p}dx > 0,
\end{equation}
implying that $\sup \Lambda(c) \leq 0$.

Now we prove that $\lambda_{1,c}, \lambda_{2,c} \in \Lambda(c)$. Take $\lambda_n \in \Lambda(c)$ such that $\lambda_n \rightarrow \lambda_{1,c}$. Then there exists $u_n \in S_c$ solving \eqref{eq1.4} with $\lambda = \lambda_n$ and $E(u_n) = m(c)$. By Lemma \ref{lem5.2}, $u_n \in H^s_{rad}$, and, up to a subsequence, $u_n \rightarrow u$ in $H^s$, where $u$ satisfies \eqref{eq1.4} with $\lambda = \lambda_{1,c}$. Note that, by the continuity of $E$ in $H^s,$ $E(u) = m(c)$ and $Q(u) = c$. Therefore, $\lambda_{1,c} \in \Lambda(c)$. Similarly, we can show that $\lambda_{2,c} \in \Lambda(c)$.

Finally, we show that $\sup \Lambda(c) < 0$. We argue by contradiction that $\sup \Lambda(c) = 0$. By the arguments above, there is a $u \in \Lambda(c)$ such that $u$ solves \eqref{eq1.4} with $\lambda = 0$. This contradicts \eqref{eq5.8} and the proof is complete.
\end{proof}

\begin{lemma}
	\label{lem5.8}  Assume that $(h)$ holds and $2 < p < 2 + (2\theta + 4s)/N$. Then any normalized ground state $u_c$ of \eqref{eq1.4} on $S_c$ has Morse index $1$.
\end{lemma}

\begin{proof}
	On the one hand, $u_c$ solves
$$
(-\Delta)^su = \lambda u + h(|x|)|u|^{p-2}u
$$
for some $\lambda < 0$. The linearized operator at $u_c$ is
$$
L_\lambda = (-\Delta)^s - \lambda - (p-1)h(|x|)|u_c|^{p-2},
$$
and
$$
\left\langle L_\lambda u_c, u_c \right\rangle = (2-p)\int_{\mathbb{R}^N}h(|x|)|u_c|^{p}dx < 0,
$$
implying that the Morse index of $u_c$ is no less than $1$ (this is a direct corollary of \cite[Theorem XIII.2]{RS}).

On the other hand, note that $S_c$ is a $C^1$ manifold with codimension $1$. Since $u_c$ is a minimizer of $E$ constrained on $S_c$, the Morse index of $u_c$ is no more than the codimension of $S_c$, i.e. no more than $1$. Hence, $u_c$ has Morse index $1$.
\end{proof}

\begin{proof}[Proof of Theorem \ref{thm2.10}] First, we prove $(i)$. By Corollary \ref{cor5.4}, $m(c)$ is strictly  decreasing on $c > 0$. Thus $m(c)$ is differentiable at almost every $c > 0$.

Next, we prove $(ii)$. Let $c_n \rightarrow c > 0$, $u_n \in S_{c_n}$, $u \in S_c$. Note that $\sqrt{\frac{c}{c_n}}u_n \in S_c$. Thus
$$
E(u) \leq E(\sqrt{\frac{c}{c_n}}u_n),
$$
and then
\begin{eqnarray}
m(c_n) - m(c) &=& E(u_n) - E(u) \nonumber \\
&\geq& E(u_n) - E(\sqrt{\frac{c}{c_n}}u_n) \nonumber \\
&=&D_uE(u_n)((1 - \sqrt{\frac{c}{c_n}})u_n) + o_n(1 - \sqrt{\frac{c}{c_n}})\nonumber \\
&=& (1 - \sqrt{\frac{c}{c_n}})D_uE(u_n)(u_n) + o_n(c_n-c) \nonumber \\
&=& 2(\sqrt{c_n} - \sqrt{c})\sqrt{c_n}\lambda_n + o_n(c_n-c),
\end{eqnarray}
where $\lambda_n$ is the Lagrange multiplier corresponding to $u_n$. Similarly, we have
$$
E(u_n) \leq E(\sqrt{\frac{c_n}{c}}u),
$$
and then
\begin{eqnarray}
m(c_n) - m(c) &=& E(u_n) - E(u) \nonumber \\
&\leq& E(\sqrt{\frac{c_n}{c}}u) - E(u) \nonumber \\
&=&D_uE(u)((\sqrt{\frac{c_n}{c}}-1)u) + o_n(\sqrt{\frac{c_n}{c}}-1) \nonumber \\
&=& (\sqrt{\frac{c_n}{c}}-1)D_uE(u)(u) + o_n(c_n-c) \nonumber \\
&=& 2(\sqrt{c_n} - \sqrt{c})\sqrt{c}\lambda + o_n(c_n-c),
\end{eqnarray}
where $\lambda$ is the Lagrange multiplier corresponding to $u$.

When $c_n > c$, we have
\begin{equation} \label{eq5.15}
\frac{2\sqrt{c_n}}{\sqrt{c_n} + \sqrt{c}}\lambda_n + o_n(1) \leq \frac{m(c_n) - m(c)}{c_n-c} \leq \frac{2\sqrt{c}}{\sqrt{c_n} + \sqrt{c}}\lambda + o_n(1).
\end{equation}
By Lemma \ref{lem5.2}, passing to a subsequence if necessary, we may assume that $u_n \rightarrow u_{\infty}$ in $H^s$ and $\lambda_n \rightarrow \lambda_{\infty}$. Using the continuity of the functional $E$ in $H^s,$ we can conclude that $u_{\infty}$ is a normalized ground state on $S_c$ and $\lambda_{\infty}$ is the Lagrange multiplier. Thus $\lambda_{\infty} \in \Lambda(c)$ and $\lambda_{\infty} \geq \lambda_{1,c}$. By Lemma \ref{lem5.4}, $\lambda_{1,c} \in \Lambda(c)$. Therefore, we can take $\lambda = \lambda_{1,c}$. Then by \eqref{eq5.15}, $\lambda_{\infty} \leq \lambda_{1,c}$. Hence, $\lambda_{\infty} = \lambda_{1,c}$ and we deduce that
$$
\lim_{n \rightarrow \infty}\frac{m(c_n) - m(c)}{c_n-c} = \lambda_{1,c},
$$
implying that $m_+'(c) = \lambda_{1,c}$.

Similarly, when $c_n < c$, we have
\begin{equation}
\frac{2\sqrt{c}}{\sqrt{c_n} + \sqrt{c}}\lambda + o_n(1) \leq \frac{m(c_n) - m(c)}{c_n-c} \leq \frac{2\sqrt{c_n}}{\sqrt{c_n} + \sqrt{c}}\lambda_n + o_n(1),
\end{equation}
and we can prove that $m_-'(c) = \lambda_{2,c}$.

It is easy to see that $m(c)$ is differentiable at $c$ if and only if $\lambda_{1,c} = \lambda_{2,c}$, i.e. $\Lambda(c)$ has only one element. The proof of $(ii)$ is complete.

Finally, we give the proof of $(iii)$. By Lemma \ref{lem5.2}, any normalized ground state is positive after multiplying $-1$ if necessary. If there exist two positive normalized ground states $u_1, u_2$ on $S_c$ for some $c > 0$. Let $\lambda(u_1), \lambda(u_2)$ be the Lagrange multiplier with respect to $u_1, u_2$ respectively. By Lemma \ref{lem5.8}, $u_1, u_2$ have Morse index $1$. Then by Theorem \ref{thm5.1}, $\lambda(u_1) \neq \lambda(u_2)$. Hence, from $(ii)$ we deduce that $m(c)$ is not differentiable at such $c$. Furthermore, by $(i)$, $m(c)$ is differentiable at almost every $c > 0$. Thus for almost every $c > 0$, \eqref{eq1.4} has a unique positive normalized ground state on $S_c$. The proof is complete.
\end{proof}

\section{Hypotheses and main results for \eqref{eq1.1} and \eqref{eq1.2}} \label{result}

\subsection{NLS with or without a potential} \label{nls}

In this subsection, we list the hypotheses and main results for \eqref{eq1.1}. We assume the following :

\begin{itemize}
	\item[$(V)$] $V(r) \in C^1(\mathbb{R}_+) \cap C([0,+\infty))$, $2V + rV'(r)$ is increasing in $r > 0$, $\lambda_1 := \inf \sigma(-\Delta + V)$ is an eigenvalue.
	\item[$(f_1)$] $f(r,t) \in C^1(\mathbb{R}_+ \times \mathbb{R}_+ ) \cap C([0,+\infty) \times [0,+\infty))$, is nonincreasing in $r > 0$ and increasing in $t > 0$, $f(r,0) = 0$, $f(r,t) = -f(r,-t)$, $(N-2)f(r,t)t \leq 2NF(r,t)+2rF_r(r,t)$ where $F(r,t) = \int_0^tf(r,s)ds$.
	\item[$(f_2)$] $(1 + \frac{4}{N})\frac{f(r,t)}{t} + \frac{2}{N}\frac{rf_r(r,t)}{t} -f_t(r,t)$ is nonincreasing in $r > 0$ and increasing in $t > 0$.
	\item[$(f_2')$] $(1 + \frac{4}{N})\frac{f(r,t)}{t} + \frac{2}{N}\frac{rf_r(r,t)}{t} -f_t(r,t)$ is nondecreasing in $r > 0$ and decreasing in $t > 0$.
	\item[$(f_3)$] There exist $p,q \in (2,2^\ast)$, where $2^* = 2N/(N-2)$ if $N \geq 3$ and $2^* = \infty$ if $N =1,2$, such that
	$$
	\lim_{t \rightarrow 0^+}\frac{f(x,t)}{t^{p-1}} = m_1(|x|),
	$$
	uniformly with $r = |x| \geq 0$,
	$$
	\lim_{t \rightarrow +\infty}\frac{f(x,t)}{t^{q-1}} = m_2(|x|),
	$$
	uniformly with $r = |x| \geq 0$, and
	$$
	\lim_{r \rightarrow +\infty}m_1(r) = m_1(\infty) \in (0,+\infty), \quad \lim_{r \rightarrow 0}m_2(r) = m_2(0) \in (0,+\infty).
	$$
\end{itemize}

Define
\begin{align*}
	\Phi_{1,\lambda}(u) & = \frac{1}{2}\int_{\mathbb{R}^N}\left(|\nabla u|^2 + V|u|^2 - \lambda |u|^2\right)dx - \int_{\mathbb{R}^N}F(|x|,u)dx, \\
E_{1}(u) & = \frac{1}{2}\int_{\mathbb{R}^N}\left(|\nabla u|^2 + V|u|^2\right)dx - \int_{\mathbb{R}^N}F(|x|,u)dx, \\
Q(u) & = \frac{1}{2}\int_{\mathbb{R}^N}|u|^2dx, \\
H^1_{V} & = \big\{u \in H^1(\mathbb{R}^N): \int_{\mathbb{R}^N}V|u|^2dx < +\infty\big\}, \quad H^1_{rad,V} = H^1_{V} \cap L^2_{rad}, \\
S_{1,c} & = \big\{u \in H^1_{V}: Q(u) = c\big\}.
\end{align*}

Using the method developed in \cite{HS}, a global branch can be established under a non-degeneracy assumption. Hence, we take the existence of such a branch as a hypothesis here. Before stating our result, we introduce two assumptions:

\begin{itemize}
	\item[$(A1)$] There exists a $C^1$ global branch
$$
u: (-\infty,\lambda_1) \rightarrow H^1_{V,rad},
$$
such that $u_\lambda = u(\lambda)$ has Morse index $1$ and solves \eqref{eq1.1}. Furthermore, $u_\lambda$ is positive, and $u_\lambda'(r) < 0$ in $r > 0$.
    \item[$(A2)$] Let $v_\lambda = \partial_\lambda u_\lambda$. Then $v_\lambda(0) < 0$.
\end{itemize}

\begin{theorem}
	\label{thm2.1}  Assume that $(A1)$, $(A2)$, $(V)$, $(f_1)$, $(f_2)$ hold. Then $\partial_\lambda\int_{\mathbb{R}^N}u_\lambda^2dx < 0$ for all $\lambda < \lambda_1$.
\end{theorem}

We also consider the case when $V \equiv 0$, in which $\lambda_1 = 0$ and $H_V^1 = H^1(\mathbb{R}^N)$. Both the mass subcritical case and the mass supercritical case will be addressed in the following:

\begin{theorem}
	\label{thm2.2}  Let $V \equiv 0$ in this Theorem.
	\begin{itemize}
		\item[$(i)$] Assume that $(A1)$, $(A2)$, $(f_1)-(f_3)$ hold with $p < 2 + 4/N$. Then $\partial_\lambda\int_{\mathbb{R}^N}u_\lambda^2dx < 0$ for all $\lambda < 0$.
		\item[$(ii)$] Assume that $(A1)$, $(A2)$, $(f_1)$, $(f_2')$, $(f_3)$ hold with $q > 2 + 4/N$. Then $\partial_\lambda\int_{\mathbb{R}^N}u_\lambda^2dx > 0$ for all $\lambda < 0$.
	\end{itemize}	 
\end{theorem}


Next, let us focus attention on the one dimensional situation, i.e. $N = 1$. We will verify $(A1)$ and $(A2)$ by assuming:
\begin{itemize}
	\item[$(H_1)$] $V \in C^1(\mathbb{R})$ is even, $\lim_{|x|\rightarrow \infty} V(x) = 0$ and $V(0) < 0$.
	\item[$(H_2)$] $f \in C^1(\mathbb{R}^2)$ with $f(x,t) = f(-x,t), f(x,t) = -f(x,-t)$ for all $(x,t) \in \mathbb{R}^2$, $f(x, 0) = 0$ for all
$x \in \mathbb{R}$, $\lim_{t \rightarrow 0}f_t(x,t) = 0$ uniformly for $x \in \mathbb{R}$  and $-f(r,t)t \leq 2F(r,t)+2rF_r(r,t)$.
    \item[$(H_3)$] $(i)$ $f_x(x,t) - V'(x)t \leq 0$ for all $t \geq 0$ and $x \geq 0$. Furthermore there exists $x_0 > 0$ such that $0 \leq x < x_0 < x'$ implies $f(x',t) - V(x')t < f(x,t) - V(x)t$ for all $t > 0$. \\
$(ii)$ $f_t(x,t)t > f(x,t) > 0$ for all $t > 0$ and $x \in \mathbb{R}$.
    \item[$(H_3')$] $(i)$ $f_x(x,t) \leq 0$ for all $t \geq 0$ and $x \geq 0$. Furthermore there exists $x_0 > 0$ such that $0 \leq x < x_0 < x'$ implies $f(x',t)t < f(x,t)$ for all $t > 0$. \\
$(ii)$ $f_t(x,t)t > f(x,t) > 0$ for all $t > 0$ and $x \in \mathbb{R}$.
    \item[$(H_4)$] There exist positive constants $\sigma$ and $A$ as well as a function $\mathcal{A} \in C^1(\mathbb{R})$ such that
$$
\lim_{t \rightarrow 0^+}\frac{f(x,t)2^\sigma}{t^{2\sigma+1}} = \mathcal{A}(x) \geq A > 0,
$$
and
$$
\lim_{t \rightarrow \infty}\lim_{x \rightarrow \infty}\frac{f(x,t)}{t} = \infty.
$$
\end{itemize}
 
\begin{theorem}
	\label{thm2.3}  Let $N = 1$. Assume that $(H_1)-(H_4)$, $(f_2)$ hold and that $2V(x) + xV'(x)$ is increasing in $x > 0$. Then there exists $u \in C^1((-\infty,\lambda_1),H^1(\mathbb{R}))$ such that for all $\lambda < \lambda_1$, $u_\lambda = u(\lambda) > 0$ solves \eqref{eq1.1} and $\partial_\lambda\int_{\mathbb{R}}u_\lambda^2dx < 0$.
\end{theorem}

\begin{theorem}
	\label{thm2.4}	Let $N = 1$ and $V \equiv 0$. Assume that $(H_2)$, $(H_3')$, $(f_3)$ hold. Then there exists $u \in C^1((-\infty,0),H^1(\mathbb{R}))$ such that $u_\lambda = u(\lambda) > 0$ solves \eqref{eq1.1}. Moreover,
	\begin{itemize}
		\item[$(i)$] if $(f_2)$ also holds with $p < 6$, then $\partial_\lambda\int_{\mathbb{R}}u_\lambda^2dx < 0$ for all $\lambda < 0$;
		\item[$(ii)$] if $(f_2')$ also holds with $q > 6$, then $\partial_\lambda\int_{\mathbb{R}}u_\lambda^2dx > 0$ for all $\lambda < 0$.
	\end{itemize}
\end{theorem}

Note that the definition in Section \ref{int} of normalized ground state does not hold for the mass supercritical case, in which $$\inf_{u \in S_{1,c}}E_1(u) = -\infty.$$ In this situation, we call $u_c$ a \emph{normalized ground state} if and only if
\begin{equation}
  E_1(u_c) = \inf\big\{E_1(u): u \in S_{1,c}, (E_1|_{S_{1,c}})'(u) = 0\big\}.
\end{equation}

\begin{theorem}[Uniqueness of the normalized ground states]
	\label{thm2.5}
\begin{itemize}
	\item[$(i)$] Under the hypotheses of Theorem \ref{thm2.1}, or of Theorem \ref{thm2.2} $(i)$, or of Theorem \ref{thm2.2} $(ii)$. We further assume that all positive solutions of \eqref{eq1.1} belong to the global branch shown by $(A1)$. Then for any $c > 0$, \eqref{eq1.1} has at most one positive normalized ground state on $S_{1,c}$.
	\item[$(ii)$] Under the hypotheses of Theorem \ref{thm2.3}, or of Theorem \ref{thm2.4} $(i)$, or of Theorem \ref{thm2.4} $(ii)$. Then for any $c > 0$, \eqref{eq1.1} has at most one positive normalized ground state on $S_{1,c}$.
\end{itemize}
\end{theorem}

\begin{remark}
	The existence of normalized ground states is not addressed in Theorem \ref{thm2.5}. Here, we provide some examples in which a normalized ground state indeed exists.
	\begin{itemize}
		\item[$(i)$] Let $V \equiv 0$ and let $f(r,t) = f(t)$ be independent of $r$. In addition to the conditions of Theorem \ref{thm2.5}, we also assume that $p, q \in (2,2+4/N)$ where $p$ and $q$ are given by $(f_3)$. By Theorem 1.1 (i) and Theorem 1.3 (i) in \cite{Shi}, \eqref{eq1.1} has a normalized ground state on $S_{1,c}$ for any $c > 0$. Hence, \eqref{eq1.1} has a unique positive normalized ground state on $S_{1,c}$ for any $c > 0$.
		\item[$(ii)$]  Let $f(r,t) = f(t)$ be independent of $r$. In addition to the conditions of Theorem \ref{thm2.5}, we also assume that $V(r) \to 0$ as $r \to \infty$, $0 \not\equiv V \leq 0$, $p, q \in (2,2+4/N)$ where $p, q$ are given by $(f_3)$, and that the following hypotheses hold:
		\begin{itemize}
			\item There exists $\delta_1 > 0$ such that $f(t)/t$ is nondecreasing in $(0, \delta_1)$.
			\item If $N \geq 5$, then $\liminf_{t \to 0} f(t)/|t|^{N/(N-2)} > 0$.
		\end{itemize}
		By \cite[Theorem A]{IM}, \eqref{eq1.1} has a normalized ground state on $S_{1,c}$ for any $c > 0$. Hence, \eqref{eq1.1} has a unique positive normalized ground state on $S_{1,c}$ for any $c > 0$.
		\item[$(iii)$] Let $V \equiv 0$ and let $f(r,t) = f(t)$ be independent of $r$. In addition to the conditions of Theorem \ref{thm2.5}, we also assume that $p, q \in (2+4/N,2^*)$ where $p$ and $q$ are given by $(f_3)$, and that the following hypotheses hold:
		\begin{itemize}
			\item $t \mapsto \tilde F(t)/|t|^{2+4/N}$ is strictly decreasing on $(-\infty, 0)$ and strictly increasing on $(0, \infty)$, where $\tilde F(t) = f(t)t-2F(t)$.
			\item When $N \geq 3$, $f(t)t < 2^*F(t)$ for all $t \in \mathbb{R}\backslash \{0\}$.
		\end{itemize} 
		By \cite[Theorem 1.1]{JL}, \eqref{eq1.1} has a normalized ground state on $S_{1,c}$ for any $c > 0$. Hence, \eqref{eq1.1} has a unique positive normalized ground state on $S_{1,c}$ for any $c > 0$.
	\end{itemize}
\end{remark}

\begin{remark}
	Let $\Psi(t,x) = e^{-i\lambda t}u_\lambda(x)$. Then $\Psi(t,x)$ is a standing wave solution of the nonlinear Schr\"{o}dinger (NLS):
	\begin{equation} \label{eq2.1}
		\begin{cases}
			i\partial_t\Psi + \Delta \Psi - V(|x|)\Psi =  -\widetilde{f}(|x|,\Psi)
			\quad \text{in } \mathbb{R}_+ \times \mathbb{R}^N, \\[1.5\jot]
			\displaystyle
			\Psi(0,x) = \Psi_0(x),
		  \end{cases}
	\end{equation}
	with $\Psi_0(x) = u_\lambda(x)$, where $\widetilde{f}(|x|,e^{i\theta}u) = e^{i\theta}f(|x|,u), u \in \mathbb{R}$. A standing wave of this kind is said to be orbitally stable in $H_V^1(\mathbb{R}^N,\mathbb{C})$ if, whenever $\|\Psi_0(x) - u_\lambda(x)\|_{H_V^1(\mathbb{R}^N,\mathbb{C})} < \delta$, the solution $\Psi(t,x)$ of \eqref{eq2.1} with initial data $\Psi_0$ exists for all $t \geq 0$, and, for any $\epsilon > 0$, there exists $\delta > 0$ such that
	$$
	\sup_{t \geq 0}\inf_{\theta \in \mathbb{R}}\|\Psi(t,x) - e^{i\theta}u_\lambda(x)\|_{H_V^1(\mathbb{R}^N,\mathbb{C})} < \epsilon,
	$$
	where
	$$
	H_V^1(\mathbb{R}^N,\mathbb{C}) = \big\{u \in H^1(\mathbb{R}^N,\mathbb{C}): \int_{\mathbb{R}^N}V|u|^2dx < +\infty\big\}.
	$$
	Note that the local well-posedness of \eqref{eq2.1} has been established in \cite{C}. Further in the subcritical case, by \cite[Corollary 6.1.2]{Caz}, problem \eqref{eq2.1} is globally well-posed if we also assume the following conditions:
	\begin{itemize}
		\item[{\bf --}] For every $K > 0$ there exists $C(K) < \infty$ such that $|f(r,t_1) - f(r,t_2)| \leq C(K)|t_1-t_2|$ for a.e. $r > 0$ and all $t_1, t_2 \in \mathbb{R}$ such that $|t_1|, |t_2| < K$.
		\item[{\bf --}] There exist $L > 0$ and $2 < d < 2  +4/N$ such that $F(r,t) \leq L(|t|^2 + |t|^d)$ for $t\in \mathbb{R}$.
	\end{itemize}
	If $u_\lambda$ is non-degenerate, using the arguments in \cite{GSS1}, we can obtain the orbital stability/instability of $\Psi(t,x) = e^{-i\lambda t}u_\lambda(x)$ immediately:
	\begin{itemize}
		\item[$(i)$] Under the hypotheses of Theorem \ref{thm2.1} or of Theorem \ref{thm2.2} $(i)$. Let $\Psi(t,x) = e^{-i\lambda t}u_\lambda(x)$ where $u_\lambda$ is on the global branch shown by $(A1)$. We further assume that $u_\lambda$ is non-degenerate in $H_V^1$. Then $\Psi$ is orbitally stable.
		\item[$(ii)$] Under the hypotheses of Theorem \ref{thm2.2} $(ii)$. Let $\Psi(t,x) = e^{-i\lambda t}u_\lambda(x)$ where $u_\lambda$ is on the global branch shown by $(A1)$. We further assume that $u_\lambda$ is non-degenerate in $H^1$. Then $\Psi$ is orbitally unstable.
	\end{itemize} 
\end{remark}

\subsection{NLS with inhomogeneous nonlinearities on the unit ball} \label{nlsb}

In this subsection, we list the hypotheses and the main results for \eqref{eq1.2}. Let
$$
L^p(B_1,|x|^{-k}) = \big\{u \in L^p(B_1): \int_{B_1}|x|^{-k}|u|^pdx < \infty\big\}.
$$
Recall that $B_1 \subset \mathbb{R}^N$ is the unit ball with $N \geq 3$, and $0 < k < 2$. When $2 \leq p \leq 2(N-k)/(N-2)$, then
$$
H_{0}^1(B_1) \hookrightarrow L^p(B_1,|x|^{-k}),
$$
and above embedding is compact provided $p < 2(N-k)/(N-2)$, see \cite[Lemma 3.2]{GY}.

Define
\begin{align*}
	\Phi_{2,\lambda}(u) & = \frac{1}{2}\int_{B_1}\left(|\nabla u|^2 - \lambda |u|^2\right)dx - \frac{1}{p}\int_{B_1}|x|^{-k}|u|^{p}dx, \\
	E_{2}(u) & = \frac{1}{2}\int_{B_1}|\nabla u|^2dx - \frac{1}{p}\int_{B_1}|x|^{-k}|u|^{p}dx, \\
Q(u) & = \frac{1}{2}\int_{B_1}|u|^2dx, \\
m_2(c) & = \inf_{u \in S_{2,c}}E_2(u) \quad \text{where } S_{2,c} = \big\{u \in H_{0}^1(B_1): Q(u) = c\big\}.
\end{align*}
For any $\lambda < \lambda_1(B_1)$, we can find a solution $u_\lambda \in H_{0}^1(B_1)$ for \eqref{eq1.2} by using the Nehari manifold method or the Mountain-pass lemma. Moreover, $u_\lambda$ is positive, radial, decreasing in $r = |x| > 0$, and has Morse index $1$. Along the lines of \cite{HS}, we can prove that
$$
\ker(-\Delta - \lambda - (p-1)|x|^{-k}|u_\lambda|^{p-2}) \cap L^2(B_1) = 0, \quad (\lambda \neq 0),
$$
and more discussions are needed on the case when $\lambda = 0$. Then using a continuation argument developed in \cite{HS}, we can show that the positive solution with Morse index $1$ is unique for any fixed $\lambda < \lambda_1(B_1)$. This is to say that for any $\lambda < \lambda_1(B_1)$, \eqref{eq1.2} has a unique, positive least action solution $u_\lambda \in H_{0,rad}^1(B_1)$, which has Morse index $1$, $u_\lambda'(r) < 0$ in $r > 0$, and $(\lambda, u_\lambda)$ is a $C^1$ curve in $\mathbb{R} \times H_{0,rad}^1(B_1)$. Recall that a solution $u$ of \eqref{eq1.2} is non-degenerate in $H^1_0(B_1)$ if for any solution of the linearized equation
\begin{equation}
- \Delta v  - \lambda v - (p-1)|x|^{-k}|u_\lambda|^{p-2}v = 0 \quad \text{in } B_1,
\end{equation}
then $v \equiv 0$.

\begin{theorem}
	\label{thm2.7} Let $N \geq 3$, $0 < k < 2$, and $2 < p < 2+2(2-k)/N$.
	Then there exists a $C^1$ global branch
$$
u: (-\infty,\lambda_1(B_1)) \rightarrow H_{0,rad}^1(B_1),
$$
such that $u_\lambda = u(\lambda)$ is positive, non-degenerate, decreasing in $r > 0$, has Morse index $1$ and solves \eqref{eq1.2}. Furthermore, all positive solutions of \eqref{eq1.2} with Morse index $1$ are on this branch, and
$\partial_\lambda\int_{B_1}|u_\lambda|^2dx < 0, \forall \lambda < \lambda_{1}(B_1).$
\end{theorem}

\begin{theorem}[Uniqueness of radial normalized ground states]
	\label{thm2.8} Under the hypotheses of Theorem \ref{thm2.7}, \eqref{eq1.2} has a unique positive normalized ground state on $S_{2,c}$ for any $c > 0$.
\end{theorem}

\begin{remark}
	The case when $p = 2 + 2(2 - k)/N$ can be addressed but it needs additional discussions. We do not provide the details here. It is worth noting that the normalized ground state exists only when $c \in (0,\hat{c})$ for some $\hat{c} < +\infty$, rather than all $c > 0$.
\end{remark}

Let $\Psi(t,x) = e^{-i\lambda t}u_\lambda(x)$. Then $\Psi(t,x)$ is a standing wave solution of the inhomogeneous nonlinear Schr\"{o}dinger on the unit ball:
\begin{equation} \label{eq2.2}
	\begin{cases}
		i\partial_t\Psi + \Delta \Psi = -|x|^{-k}|\Psi|^{p-2}\Psi
		\quad \text{in } \mathbb{R}_+ \times B_1, \\[1.5\jot]
		\displaystyle
		\Psi(0,x) = \Psi_0(x),
	  \end{cases}
\end{equation}
with $\Psi_0(x) = u_\lambda(x)$. A standing wave of this kind is said to be orbitally stable in $H_0^1(B_1,\mathbb{C})$ if, whenever $\|\Psi_0(x) - u_\lambda(x)\|_{H_0^1(B_1,\mathbb{C})} < \delta$, the solution $\Psi(t,x)$ of \eqref{eq2.2} with initial data $\Psi_0$ exists for all $t \geq 0$, and, for any $\epsilon > 0$, there exists $\delta > 0$ such that
$$
\sup_{t \geq 0}\inf_{\theta \in \mathbb{R}}\|\Psi(t,x) - e^{i\theta}u_\lambda(x)\|_{H_0^1(B_1,\mathbb{C})} < \epsilon.
$$
Finally in this subsection, we state our result of orbital stability of standing wave solutions.

\begin{theorem}[Orbital stability]
	\label{thm2.9} Under the hypotheses of Theorem \ref{thm2.7}. Let $\Psi(t,x) = e^{-i\lambda t}u_\lambda(x)$ where $u_\lambda$ is on the global branch given by Theorem \ref{thm2.7}. We further assume that \eqref{eq2.2} is locally well-posed. Then $\Psi$ is orbitally stable.
\end{theorem}

\begin{remark}
	Theorem \ref{thm2.9} generalizes the result in \cite[Theorem 1.7 (1)]{NTV} where $k = 0$ was considered.
\end{remark}

\section{Proofs of Theorems \ref{thm2.1}-\ref{thm2.5}} \label{pr1}

\subsection{Monotonicity of $L^2$ mass of solutions on the global branch for \eqref{eq1.1}}

Let $v_\lambda = \partial_\lambda u_\lambda$, where $u_\lambda$ is given by Subsection \ref{nls}. Then $v_\lambda$ is radially symmetric, satisfying
\begin{equation} \label{eqB.5}
  -\Delta v_\lambda + Vv_\lambda = \lambda v_\lambda + u_\lambda + f_t(|x|,u_\lambda)v_\lambda.
\end{equation}

\begin{lemma}
\label{lemB.1}  Assume that $(A1)$ and $(A2)$ hold. Then $v_\lambda = v_\lambda(r)$ changes sign at most once in $r > 0$.
\end{lemma}

\begin{proof}
	We argue by contradiction. By $(A2)$, we may assume that $v_\lambda < 0$ in $(0,r_1)$ and $(r_2,r_3)$ for some $0 < r_1 < r_2 < r_3 \leq +\infty$. Set
$$
v_1 = \chi_{(0,r_1)}v_\lambda, \quad v_2 = \chi_{(r_2,r_3)}v_\lambda.
$$
Note that $v_1$ and $v_2$ are linear independent. By (\ref{eqB.5}), we have
$$
\int_{\mathbb{R}^N}(|\nabla v_1|^2 + Vv_1^2 - \lambda v_1^2 - f_t(|x|,u_\lambda)v_1^2)dx = \int_{\mathbb{R}^N}u_\lambda v_1dx < 0,
$$
$$
\int_{\mathbb{R}^N}(|\nabla v_2|^2 + Vv_2^2 - \lambda v_2^2 - f_t(|x|,u_\lambda)v_2^2)dx = \int_{\mathbb{R}^N}u_\lambda v_2dx < 0,
$$
which is in contradiction with the fact that the Morse index of $u_\lambda$ is $1$. Hence, we arrive at our conclusion.
\end{proof} 

\begin{proof}[Proof of Theorem \ref{thm2.1}]  

The proof is divided into two steps.

\vskip0.1in
\emph{Step 1:} $\int_{\mathbb{R}^N}u_\lambda v_\lambda dx \neq 0$ for any $\lambda < \lambda_1$.

Since $u_\lambda$ solves (\ref{eq1.1}), we have
\begin{equation} \label{eqB.1}
\int_{\mathbb{R}^N}(\nabla u_\lambda \nabla v_\lambda + Vu_\lambda v_\lambda) dx = \lambda\int_{\mathbb{R}^N}u_\lambda v_\lambda dx + \int_{\mathbb{R}^N}f(|x|,u_\lambda)v_\lambda dx,
\end{equation}
and the following Pohozaev identity holds
\begin{eqnarray} \label{eqB.2}
&& \frac{N-2}{2}\int_{\mathbb{R}^N}|\nabla u_\lambda|^2 dx + \frac{N}{2}\int_{\mathbb{R}^N}Vu_\lambda^2 dx + \frac{1}{2}\int_{\mathbb{R}^N}|x|V'(|x|) u_\lambda^2 dx \nonumber \\
&=& \frac{\lambda N}{2}\int_{\mathbb{R}^N}u_\lambda^2 dx + N\int_{\mathbb{R}^N}F(|x|,u_\lambda)dx + \int_{\mathbb{R}^N}|x|F_r(|x|,u_\lambda)dx.
\end{eqnarray}
Differentiating both sides of (\ref{eqB.2}) with respect to $\lambda$, we have
\begin{eqnarray} \label{eqB.3}
&& (N-2)\int_{\mathbb{R}^N}\nabla u_\lambda \nabla v_\lambda dx + \int_{\mathbb{R}^N}(NV +|x|V'(|x|))u_\lambda v_\lambda dx = \frac{N}{2}\int_{\mathbb{R}^N}u_\lambda^2dx \nonumber \\
&+& \lambda N\int_{\mathbb{R}^N}u_\lambda v_\lambda dx + N\int_{\mathbb{R}^N}f(|x|,u_\lambda)v_\lambda dx + \int_{\mathbb{R}^N}|x|f_r(|x|,u_\lambda)v_\lambda dx.
\end{eqnarray}
The combination of (\ref{eqB.1}) and (\ref{eqB.3}) imply that
\begin{eqnarray} \label{eqB.4}
&& \int_{\mathbb{R}^N}(2V +|x|V'(|x|))u_\lambda v_\lambda dx - 2\lambda \int_{\mathbb{R}^N}u_\lambda v_\lambda dx \nonumber \\
&=&  \frac{N}{2}\int_{\mathbb{R}^N}u_\lambda^2dx + 2\int_{\mathbb{R}^N}f(|x|,u_\lambda)v_\lambda dx + \int_{\mathbb{R}^N}|x|f_r(|x|,u_\lambda)v_\lambda dx.
\end{eqnarray}
By (\ref{eqB.5}),
\begin{equation} \label{eqB.6}
\int_{\mathbb{R}^N}(\nabla u_\lambda \nabla v_\lambda + Vu_\lambda v_\lambda) dx = \lambda \int_{\mathbb{R}^N} u_\lambda v_\lambda dx + \int_{\mathbb{R}^N}u_\lambda^2 dx + \int_{\mathbb{R}^N}f_t(|x|,u_\lambda)u_\lambda v_\lambda dx.
\end{equation}
By (\ref{eqB.1}) and (\ref{eqB.6}),
\begin{equation} \label{eqB.7}
\int_{\mathbb{R}^N}u_\lambda^2 dx = \int_{\mathbb{R}^N}f(|x|,u_\lambda)v_\lambda dx - \int_{\mathbb{R}^N}f_t(|x|,u_\lambda)u_\lambda v_\lambda dx.
\end{equation}
Then from (\ref{eqB.4}) and (\ref{eqB.7}), we derive that
\begin{eqnarray} \label{eqB.8}
&& \int_{\mathbb{R}^N}\left(\frac{N+4}{2}\frac{f(|x|,u_\lambda)}{u_\lambda} + \frac{|x|f_r(|x|,u_\lambda)}{u_\lambda} - \frac{N}{2}f_t(|x|,u_\lambda) - \left(2V +|x|V'(|x|)\right)\right)u_\lambda v_\lambda dx \nonumber \\
&=& -2\lambda \int_{\mathbb{R}^N}u_\lambda v_\lambda dx.
\end{eqnarray}

Let us argue by contradiction, and assume that $\int_{\mathbb{R}^N}u_\lambda v_\lambda dx = 0$. This implies that $v_\lambda$ is sign-changing. By Lemma \ref{lemB.1}, we know that $v_\lambda = v_\lambda(r)$ changes sign exactly once in $r > 0$. We assume that $v_\lambda(r) < 0$ in $(0,r^\ast)$ and $v_\lambda(r) \geq 0$ in $(r^\ast,+\infty)$.

From (\ref{eqB.8}), $(V)$ and $(f_2)$, we derive a self-contradictory inequality
\begin{eqnarray}
0 &=& -2\lambda \int_{\mathbb{R}^N}u_\lambda v_\lambda dx \nonumber \\
&=& \int_{\mathbb{R}^N}\left(\frac{N+4}{2}\frac{f(|x|,u_\lambda)}{u_\lambda} + \frac{|x|f_r(|x|,u_\lambda)}{u_\lambda} - \frac{N}{2}f_t(|x|,u_\lambda) - \left(2V +|x|V'(|x|)\right)\right)u_\lambda v_\lambda dx \nonumber \\
&<& C(r^\ast)\int_{B_{r^\ast}}u_\lambda v_\lambda dx + C(r^\ast)\int_{B_{r^\ast}^c}u_\lambda v_\lambda dx \nonumber \\
&=& C(r^\ast)\int_{\mathbb{R}^N}u_\lambda v_\lambda dx = 0,
\end{eqnarray}
where
$$
C(r^\ast) = \frac{N+4}{2}\frac{f(r^\ast,u_\lambda(r^\ast))}{u_\lambda(r^\ast)} + \frac{|x|f_r(r^\ast,u_\lambda(r^\ast))}{u_\lambda(r^\ast)} - \frac{N}{2}f_t(r^\ast,u_\lambda(r^\ast)) - (2V(r^\ast) +r^\ast V'(r^\ast)).
$$

\vskip0.1in
\emph{Step 2:} Completion of the proof.

Noticing that $\partial_\lambda\int_{\mathbb{R}^N}u_\lambda^2dx$ is continuous with respect to $\lambda < \lambda_1$ by $(A1)$, we have $\partial_\lambda\int_{\mathbb{R}^N}u_\lambda^2dx < 0$, or $\partial_\lambda\int_{\mathbb{R}^N}u_\lambda^2dx > 0$ for all $\lambda < \lambda_1$ by Step 1. Then classical bifurcation arguments (c.f. \cite[Theorem 1.7]{CR}) yield that
$$
\lim_{\lambda \rightarrow \lambda_1}\int_{\mathbb{R}^N}u_\lambda^2dx = 0
$$
implying that $\partial_\lambda\int_{\mathbb{R}^N}u_\lambda^2dx < 0$. The proof is complete.
\end{proof}

\begin{proof}[Proof of Theorem \ref{thm2.2}]  

The proof is divided into two steps.

\vskip0.1in
\emph{Step 1:} $\int_{\mathbb{R}^N}u_\lambda v_\lambda dx \neq 0$ for any $\lambda < 0$.

Let $V = 0$ in \eqref{eqB.8}, we derive that
\begin{equation} \label{eqA.11}
-2\lambda \int_{\mathbb{R}^N}u_\lambda v_\lambda dx =  \int_{\mathbb{R}^N}\left(\frac{N+4}{2}\frac{f(|x|,u_\lambda)}{u_\lambda} + \frac{|x|f_r(|x|,u_\lambda)}{u_\lambda} - \frac{N}{2}f_t(|x|,u_\lambda)\right)u_\lambda v_\lambda dx.
\end{equation}
We will argue by contradiction, and assume that $\int_{\mathbb{R}^N}u_\lambda v_\lambda dx = 0$. This implies that $v_\lambda$ is sign-changing and lead to the contradiction. First we assume that $(f_2)$ holds and prove $(i)$. Thanks to Lemma \ref{lemB.1}, we know that $v_\lambda = v_\lambda(r)$ changes sign exactly once in $r > 0$. We assume that $v_\lambda(r) < 0$ in $(0,r^\ast)$ and $v_\lambda(r) \geq 0$ in $(r^\ast,+\infty)$.

If $(f_2)$ holds, from (\ref{eqA.11}), we derive a self-contradictory inequality
\begin{eqnarray}
0 &=& \int_{\mathbb{R}^N}\left(\frac{N+4}{2}\frac{f(|x|,u_\lambda)}{u_\lambda} + \frac{|x|f_r(|x|,u_\lambda)}{u_\lambda} - \frac{N}{2}f_t(|x|,u_\lambda)\right)u_\lambda v_\lambda dx \nonumber \\
&<& C(r^\ast)\int_{B_{r^\ast}}u_\lambda v_\lambda dx + C(r^\ast)\int_{B_{r^\ast}^c}u_\lambda v_\lambda dx \nonumber \\
&=& C(r^\ast)\int_{\mathbb{R}^N}u_\lambda v_\lambda dx = 0,
\end{eqnarray}
where
$$
C(r^\ast) = \frac{N+4}{2}\frac{f(r^\ast,u_\lambda(r^\ast))}{u_\lambda(r^\ast)} + \frac{|x|f_r(r^\ast,u_\lambda(r^\ast))}{u_\lambda(r^\ast)} - \frac{N}{2}f_t(r^\ast,u_\lambda(r^\ast)).
$$

Similarly, if $(f_2')$ holds, from (\ref{eqA.11}), we derive a self-contradictory inequality
\begin{eqnarray}
0 &=& \int_{\mathbb{R}^N}\left(\frac{N+4}{2}\frac{f(|x|,u_\lambda)}{u_\lambda} + \frac{|x|f_r(|x|,u_\lambda)}{u_\lambda} - \frac{N}{2}f_t(|x|,u_\lambda)\right)u_\lambda v_\lambda dx \nonumber \\
&>& C(r^\ast)\int_{B_{r^\ast}}u_\lambda v_\lambda dx + C(r^\ast)\int_{B_{r^\ast}^c}u_\lambda v_\lambda dx \nonumber \\
&=& C(r^\ast)\int_{\mathbb{R}^N}u_\lambda v_\lambda dx = 0.
\end{eqnarray}

\vskip0.1in
\emph{Step 2:} Completion of the proof.

Noticing that $\partial_\lambda\int_{\mathbb{R}^N}u_\lambda^2dx$ is continuous with respect to $\lambda < 0$, we have $\partial_\lambda\int_{\mathbb{R}^N}u_\lambda^2dx < 0$ , or $\partial_\lambda\int_{\mathbb{R}^N}u_\lambda^2dx > 0$ for all $\lambda < 0$ by Step 1. If $(f_3)$ holds with $p < 2 + 4/N$, \cite[Corollary 4.2]{HS} yields that
$$
\lim_{\lambda \rightarrow 0^-}\int_{\mathbb{R}^N}u_\lambda^2dx = 0,
$$
implying that $\partial_\lambda\int_{\mathbb{R}^N}u_\lambda^2dx < 0$. Similarly, If $(f_3)$ holds with $q > 2 + 4/N$, \cite[Corollary 4.5]{HS} yields that
$$
\lim_{\lambda \rightarrow -\infty}\int_{\mathbb{R}^N}u_\lambda^2dx = 0,
$$
implying that $\partial_\lambda\int_{\mathbb{R}^N}u_\lambda^2dx > 0$. The proof is complete.
\end{proof}

\begin{proof}[Proof of Theorem \ref{thm2.3}]
When $N = 1$, under the hypotheses $(H_1)-(H_4)$, \cite{JS} showed that $(A1)$ and $(A2)$ hold. Obviously, $(V)$ and $(f_1)$ hold. Then along the lines of the proof of Theorem \ref{thm2.1}, we can complete the proof of Theorem \ref{thm2.3}.
\end{proof}

\begin{remark}
	A similar result to Theorem \ref{thm2.3} has been shown in \cite{MST}.
\end{remark}

\begin{proof}[Proof of Theorem \ref{thm2.4}]
	Mimicking the proofs in \cite{JS}, we know that $(A1)$ and $(A2)$ hold under the hypotheses of Theorem \ref{thm2.4}. Along the lines of the proof of Theorem \ref{thm2.2}, we can complete the proof of Theorem \ref{thm2.4}.
\end{proof}

\subsection{Uniqueness of the normalized ground state for \eqref{eq1.1}}

\begin{proof}[Proof of Theorem \ref{thm2.5}]  

	The proof is divided into three steps.

	\vskip0.1in
\emph{Step 1:} Any normalized ground state of \eqref{eq1.1} is positive after multiplying $-1$ if necessary.

If $u_c$ satisfies $E_1(u_c) = \inf_{u \in S_{1,c}}E_1(u)$. Noticing that $E_1(|u_c|) \leq E_1(u_c)$, we know that $|u_c|$ is also a normalized ground state, which solves \eqref{eq1.1} for some $\lambda$. By the strong maximum principle, $|u_c| > 0$. Since $u_c$ is continuous, we know $u_c$ cannot change its sign. After multiplying $-1$ if necessary, we may assume that $u_c = |u_c| > 0$.

\vskip0.1in
\emph{Step 2:} If $u > 0$ solves \eqref{eq1.1} for some $\lambda$, then $\lambda < \lambda_1$.

First, we consider the case when $(V)$ holds. Let $e_1$ be the positive, unit eigenfunction corresponding to $\lambda_1$. Then we have
\begin{equation}
\int_{\mathbb{R}^N}f(|x|,u)e_1dx = \int_{\mathbb{R}^N}(\nabla u \nabla e_1 + Vu e_1 - \lambda u e_1)dx = (\lambda_1 - \lambda)\int_{\mathbb{R}^N} u e_1dx.
\end{equation}
By $(f_1)$, $f(|x|,u) > 0$ when $u > 0$. Hence, $\int_{\mathbb{R}^N}f(|x|,u)e_1dx > 0$. Noticing that $\int_{\mathbb{R}^N} u e_1dx > 0$, we know that $\lambda < \lambda_1$.

Next, let $V \equiv 0$ and $\lambda_1 = 0$. The Pohozaev identity yields that
\begin{equation} \label{eq3.14}
2\lambda\int_{\mathbb{R}^N} |u|^2dx = \int_{\mathbb{R}^N}\left( (N-2)f(|x|,u)u - 2NF (|x|,u)+ 2|x|F_r(|x|,u)\right) dx.
\end{equation}
By $(f_1)$, the right part of \eqref{eq3.14} is less than $0$. Therefore, $\lambda < 0$.

\vskip0.1in
\emph{Step 3:} Completion of the proof.

Suppose on the contrary that there are two positive normalized ground states $u_{c,1}, u_{c,2}$ with same $L^2$ mass $c$. Then $u_{c,1}, u_{c,2}$ solve \eqref{eq1.1} for some $\lambda, \widehat{\lambda}$ respectively. By Step 2, $\lambda < \lambda_1, \widehat{\lambda} < \lambda_1$. Thus both $u_{c,1}, u_{c,2}$ are on the global branch. Then, the fact $\int_{\mathbb{R}^N}|u_{c,1}|^2dx = \int_{\mathbb{R}^N}|u_{c,2}|^2dx$ implies that $\lambda = \widehat{\lambda}$ since the $L^2$ norm of solutions on the global branch is monotonic. Then we derive that $u_{c,1} = u_{c,2}$, completing the proof of Theorem \ref{thm2.5}.
\end{proof}



\section{Proofs of Theorems \ref{thm2.7}-\ref{thm2.9}} \label{pr2}

\subsection{Existence and uniqueness of the global branch and monotonicity of $L^2$ mass of solutions on this branch for \eqref{eq1.2}}

\begin{lemma}\label{lemma 5.1}
	\label{lem4.1} Let $N \geq 3$, $0 < k < 2$, and $2 < p < 2+2(2-k)/N$. If $u_\lambda \in H_{0,rad}^1(B_1)$ is a positive solution of \eqref{eq1.2} with $\lambda < \lambda(B_1)$, which has Morse index $1$ and decreasing in $r > 0$. Then $u_\lambda$ is non-degenerate in $H_0^1(B_1)$.
\end{lemma}

\begin{proof} 
	Let $w \in H^1_0(B_1)$ satisfy
\begin{equation} \label{eq4.1}
- \Delta w  - \lambda w - (p-1)|x|^{-k}|u_\lambda|^{p-2}w = 0 \ in \ B_1.
\end{equation}
Since $u_\lambda \in H^1_{0,rad}(B_1)$, we know that $w \in H^1_{0,rad}(B_1)$. If $w \neq 0$, since the Morse index of $u_\lambda$ is $1$, $0$ is the second eigenvalue of $- \Delta  - \lambda - (p-1)|x|^{-k}|u_\lambda|^{p-2}$ with the form domain $H^1_{0,rad}(B_1)$ and $w$ is the corresponding eigenfunction. Hence, $w$ changes sign exactly once. Direct calculation yields to
\begin{equation} \label{eq4.2}
- \Delta u_\lambda  - \lambda u_\lambda - (p-1)|x|^{-k}|u_\lambda|^{p-2}u_\lambda = (2-p)|x|^{-k}|u_\lambda|^{p-2}u_\lambda,
\end{equation}
\begin{equation} \label{eq4.3}
- \Delta (x \cdot \nabla u_\lambda)  - \lambda x \cdot \nabla u_\lambda - (p-1)|x|^{-k}|u_\lambda|^{p-2}x \cdot \nabla u_\lambda = (2-k)|x|^{-k}|u_\lambda|^{p-2}u_\lambda + 2\lambda u_\lambda.
\end{equation}
Next two cases will be treated.

\vskip0.1in
\emph{Case 1:} $\lambda \neq 0$.

Since $w$ changes sign exactly once, we may assume that $w(r_0) = 0$, $w(r) > 0$ in $r < r_0$ and $w(r) \leq 0$ in $r_0 < r < 1$ for some $r_0 \in (0,1)$. Set
$$
\phi(x) = \lambda u_\lambda + m|x|^{-k}|u_\lambda|^{p-2}u_\lambda
$$
where $m$ is chosen such that $\phi(x) = 0$ for $|x| = r_0$. On the one hand, \eqref{eq4.2} and \eqref{eq4.3} imply that $\phi$ is in the range of the operator $- \Delta  - \lambda - (p-1)|x|^{-k}|u_\lambda|^{p-2}$ and thus $\int_{B_1}\phi wdx = 0$. On the other hand, $\phi/u_\lambda$ is monotonic in $r \in (0,1)$. Hence, $(\phi/u_\lambda) w$ does not change sign in $r \in (0,1)$, so that$\int_{B_1}\phi wdx > 0$, which is a contradiction. Therefore, $w = 0$.

\vskip0.1in
\emph{Case 2:} $\lambda = 0$.

Note that $|x|^{-k}|u_\lambda|^{p-2}u_\lambda$ is in the range of the operator $- \Delta - (p-1)|x|^{-k}|u_\lambda|^{p-2}$. Then $\int_{B_1}|x|^{-k}|u_\lambda|^{p-2}u_\lambda wdx = 0$. Applying Green's Theorem and \eqref{eq4.2}, \eqref{eq4.3}, we have
\begin{eqnarray} \label{eq4.4}
\int_{\partial B_1}x \cdot \nabla u_\lambda \frac{\partial w}{\partial n} dS &=& \int_{B_1}(x \cdot \nabla u_\lambda \Delta w - \Delta (x \cdot \nabla u_\lambda) w)dx \nonumber \\
&=& (2-k)\int_{B_1}|x|^{-k}|u_\lambda|^{p-2}u_\lambda wdx \nonumber \\
&=& 0.
\end{eqnarray}
On $\partial B_1$, $x \cdot \nabla u_\lambda = u_\lambda'(1) \neq 0$. Thus $\frac{\partial w}{\partial n} = w'(1) = 0$, which is impossible by Hopf's lemma. Hence, $w = 0$.
\end{proof}

\begin{proof}[Proof of Theorem \ref{thm2.7}]
	Thanks to Lemma \ref{lem4.1}, by Theorem \ref{exi}, we can obtain the existence of a $C^1$ global branch
$$
u: (-\infty,\lambda_1(B_1)) \rightarrow H_{0,rad}^1(B_1),
$$
such that $u_\lambda = u(\lambda)$ is positive, non-degenerate, has Morse index $1$ and solves \eqref{eq1.2} (see more details in the proof of \cite[Theorem 5.11]{HS}). Using the moving plane method, we can show that $u_\lambda$ is decreasing in $r \in (0,1)$ (see e.g. \cite{GNN}). Along the lines of \cite[Theorem 3.2]{HS}, we can prove that $\int_{B_1}u_\lambda^2 dx \rightarrow \infty$ as $\lambda \rightarrow -\infty$. Noticing that $\partial_\lambda\int_{B_1}u_\lambda^2dx$ is continuous with respect to $\lambda < \lambda_1$, we have $\partial_\lambda\int_{B_1}u_\lambda^2dx < 0$ for all $\lambda < \lambda_1(B_1)$ by the following Lemma \ref{lemD.2}. In the remaining  part of this subsection, we further need to show the uniqueness of such a global branch to complete the proof of Theorem \ref{thm2.7}. 
\end{proof}

Let $v_\lambda = \partial_\lambda u_\lambda$. Then $v_\lambda \in H_0^1(B_1)$ is radially symmetric, satisfying
\begin{equation} \label{eqD.6}
	-\Delta v_\lambda = \lambda v_\lambda + u_\lambda + (p-1)|x|^{-k}|u_\lambda|^{p-2}v_\lambda.
\end{equation}
We do not aim to show that $v_\lambda(0) < 0$ like $(A2)$ assumed in Theorems \ref{thm2.1}, \ref{thm2.2}. Instead, we will prove that $v_\lambda(r) < 0$ near $r = 1$ if $\int_{B_1}u_\lambda v_\lambda dx = 0$.

\begin{lemma}
	\label{lemD.2} Let $N \geq 3$, $0 < k < 2$, and $2 < p < 2+2(2-k)/N$. Then $\int_{B_1}u_\lambda v_\lambda dx \neq 0$ for all $\lambda < \lambda_{1}(B_1)$.
\end{lemma}

\begin{proof} 
	Suppose on the contrary that $\int_{B_1}u_\lambda v_\lambda dx = 0$ for some $\lambda < \lambda_{1}(B_1)$. Since $u_\lambda$ solves (\ref{eq1.2}), we have
\begin{equation} \label{eqD.2}
\int_{B_1}\nabla u_\lambda \nabla v_\lambda dx  =  \int_{B_1}|x|^{-k}|u_\lambda|^{p-2}u_\lambda v_\lambda dx,
\end{equation}
and the following Pohozaev identity holds
\begin{equation} \label{eqD.3}
\frac{N-2}{2}\int_{B_1}|\nabla u_\lambda|^2 dx + \frac{1}{2}\int_{\partial B_1}|\frac{\partial u_\lambda}{\partial n}|^2 dS = \frac{\lambda N}{2}\int_{B_1}|u_\lambda|^2 dx + \frac{N-k}{p}\int_{B_1}|x|^{-k}|u_\lambda|^pdx.
\end{equation}
Differentiating both sides of (\ref{eqD.3}) with respect to $\lambda$, we have
\begin{eqnarray} \label{eqD.4}
&& (N-2)\int_{B_1}\nabla u_\lambda \nabla v_\lambda dx + \omega_Nu'_\lambda(1)v'_\lambda(1) \nonumber \\
&=& \frac{N}{2}\int_{B_1}|u_\lambda|^2dx + (N-k)\int_{B_1}|x|^{-k}|u_\lambda|^{p-2}u_\lambda v_\lambda dx.
\end{eqnarray}
The combination of (\ref{eqD.2}) and (\ref{eqD.4}) imply that
\begin{eqnarray} \label{eqD.5}
\omega_Nu'_\lambda(1)v'_\lambda(1) =  \frac{N}{2}\int_{B_1}|u_\lambda|^2dx + (2-k)\int_{B_1}|x|^{-k}|u_\lambda|^{p-2}u_\lambda v_\lambda.
\end{eqnarray}
Since $v_\lambda$ satisfies (\ref{eqD.6}), we have
\begin{equation} \label{eqD.7}
\int_{B_1}\nabla u_\lambda \nabla v_\lambda dx  =  \int_{B_1}|u_\lambda|^2dx + (p-1)\int_{B_1}|x|^{-k}|u_\lambda|^{p-2}u_\lambda v_\lambda dx.
\end{equation}
From (\ref{eqD.2}), (\ref{eqD.5}) and (\ref{eqD.7}), we derive that
\begin{eqnarray} \label{eqD.8}
\omega_Nu'_\lambda(1)v'_\lambda(1) =  (\frac{N}{2} - \frac{2-k}{p-2})\int_{B_1}|u_\lambda|^2dx < 0.
\end{eqnarray}
By the strong maximum principle, $u'_\lambda(1) > 0$. Hence, $v'_\lambda(1) < 0$, implying that $v_\lambda = v_\lambda(r) < 0$ in $(1-\epsilon,1)$ for some $\epsilon > 0$. Then similar to the proof of Lemma \ref{lemB.1}, we can prove that $v_\lambda(r) \leq 0$ in $(0,1)$ or $v_\lambda(r) \geq 0$ in $(0,r_0)$ and $v_\lambda(r) < 0$ in $(r_0,1)$ for some $r_0 \in (0,1)$. On the one hand,
\begin{eqnarray} \label{eqD.9}
&& \int_{B_1}|x|^{-k}|u_\lambda|^{p-2}u_\lambda v_\lambda dx \nonumber \\
&>& r_0^{-k}|u_\lambda(r_0)|^{p-2}\int_{B_{r_0}}u_\lambda v_\lambda dx + r_0^{-k}|u_\lambda(r_0)|^{p-2}\int_{r_0 < |x| < 1}u_\lambda v_\lambda dx \nonumber \\
&=& 0.
\end{eqnarray}
On the other hand, by (\ref{eqD.2}) and (\ref{eqD.7}),
\begin{equation} \label{eqD.10}
(2-p)\int_{B_1}|x|^{-k}|u_\lambda|^{p-2}u_\lambda v_\lambda dx  =  \int_{B_1}|u_\lambda|^2dx > 0,
\end{equation}
in a contradiction with (\ref{eqD.9}). The proof is complete.
\end{proof}

\begin{proof}[Completion of the proof for Theorem \ref{thm2.7}]
	Arguing by contradiction, we assume that there exists $\widehat{u}_{\lambda_0} \in H^1_{0}(B_1)$, which is positive, has Morse index $1$ and solves \eqref{eq1.2} with $\lambda = \lambda_0 < \lambda_1(B_1)$, but is not on the global branch $u(\lambda)$. Using the moving plane method, we can show that $\widehat{u}_{\lambda_0}$ is radial and decreasing in $r \in (0,1)$. Then by Lemma \ref{lem4.1}, $\widehat{u}_{\lambda_0}$ is non-degenerate. Along the lines of \cite[Lemma 2.5]{HS}, there exists
$$
\widehat{u}: (-\infty,\lambda_0] \rightarrow H_{0,rad}^1(B_1),
$$
such that $\widehat{u}_{\lambda_0} = \widehat{u}(\lambda_0)$, $\widehat{u}(\lambda)$ is positive, non-degenerate, has Morse index $1$ and solves \eqref{eq1.2}. We assume that the branch $\widehat{u}(\lambda)$ can be extended to $(-\infty,\lambda_\ast)$.

By Lemma \ref{lemD.2}, we have $\partial_\lambda\int_{B_1}|\widehat{u}(\lambda)|^2dx < 0$ for all $\lambda < \lambda_\ast$. Note that
\begin{eqnarray}
\partial_\lambda\int_{B_1}|x|^{-k}|\widehat{u}(\lambda)|^pdx &=& p\int_{B_1}|x|^{-k}|\widehat{u}(\lambda)|^{p-2}\widehat{u}(\lambda)\partial_\lambda\widehat{u}(\lambda)dx \nonumber \\
&=& \frac{p}{2-p}\int_{B_1}|\widehat{u}(\lambda)|^2dx \nonumber \\
&\leq& 0.
\end{eqnarray}
Hence, $\int_{B_1}|\nabla \widehat{u}(\lambda)|^2dx = \lambda\int_{B_1}|\widehat{u}(\lambda)|^2dx + \int_{B_1}|x|^{-k}|\widehat{u}(\lambda)|^pdx$ is uniformly bounded when $\lambda \rightarrow \lambda_\ast^-$. Then mimicking the proof of \cite[Lemma 2.5]{HS}, we can prove that $\lambda_\ast = \lambda_1(B_1)$.

Now we aim to show that $\widehat{u}(\lambda) \rightarrow 0$ as $\lambda \to \lambda_1(B_1)^-$.
Let $\lambda_n \rightarrow \lambda_1(B_1)^-$. Noticing that $\widehat{u}(\lambda_n)$ is bounded in $H_0^1(B_1)$, up to a subsequence, we may assume that $\widehat{u}(\lambda_n) \rightharpoonup \widehat{u}$ in $H_0^1(B_1)$, $\widehat{u}(\lambda_n) \rightarrow \widehat{u}$ in $L^2(B_1)$ and $L^p(B_1,|x|^{-k})$. Note that $\widehat{u}$ is a nonnegative solution of \eqref{eq1.2} with $\lambda = \lambda_1(B_1)$. By the Step 3 in the proof of Theorem \ref{thm2.8} below, we know that $\widehat{u} = 0$. Hence, $\widehat{u}(\lambda_n) \rightarrow 0$ in $L^2(B_1)$ and $L^p(B_1,|x|^{-k})$, and then
$$
\int_{B_1}|\nabla \widehat{u}(\lambda_n)|^2dx = \lambda_n\int_{B_1}|\widehat{u}(\lambda_n)|^2dx + \int_{B_1}|x|^{-k}|\widehat{u}(\lambda_n)|^pdx \rightarrow 0.
$$
The global branch $u(\lambda)$ also stems form $\lambda_1(B_1)$. However, by the classical bifurcation theory (c.f. \cite[Theorem 1.7]{CR}), there exists only one local branch of positive, radial solutions of \eqref{eq1.2} stemming form $\lambda_1(B_1)$, which is a contradiction. The proof of Theorem \ref{thm2.7} is now complete.
\end{proof}

\subsection{Uniqueness of the normalized ground state for \eqref{eq1.2}}

\begin{proof}[Proof of Theorem \ref{thm2.8}]
	 The proof is divided into four steps.
	
	\vskip0.1in
	\emph{Step 1:} Existence of a normalized ground state.
	
	Modifying the proof of \cite[Lemma 3.2]{GY}, we have
	\begin{align*}
		\int_{B_1}|x|^{-k}|u|^pdx & = \int_{B_1}|x|^{-k}|u|^k|u|^{p-k}dx \\
		& \leq \left( \int_{B_1}|x|^{-2}|u|^2dx\right)^{\frac k2} \left( \int_{B_1}|u|^{\frac{2(p-k)}{2-k}}dx\right)^{\frac {2-k}2} \\
		& \leq C\left( \int_{B_1}|\nabla u|^2dx\right)^{\frac k2} \left( \int_{B_1}|u|^{\frac{2(p-k)}{2-k}}dx\right)^{\frac {2-k}2} \\
		& \leq C_1\left( \int_{B_1}|\nabla u|^{2}dx\right)^{\frac k2 + \frac{N(p-2)}{4}}\left( \int_{B_1}| u|^{2}dx\right)^{\frac {p-k}2 - \frac{N(p-2)}{4}},
	\end{align*}
	where $C, C_1$ are positive constants and we have used Gagliardo-Nirenberg inequality in the final step. Note that $0 < k < 2$ and $2 < p < 2 + 2(2-k)/N$. Let $\{u_n\} \subset S_{2,c}$ be a minimizing sequence for $m_2(c)$. From
	\begin{align*}
		E_2(u_n) \geq \frac12\int_{B_1}|\nabla u_n|^{2}dx - \frac{C_1}{p}\left( \int_{B_1}|\nabla u_n|^{2}dx\right)^{\frac k2 + \frac{N(p-2)}{4}}\left( \int_{B_1}| u_n|^{2}dx\right)^{\frac {p-k}2 - \frac{N(p-2)}{4}},
	\end{align*}
    we get $\{u_n\}$ is bounded in $H_0^1(B_1)$ since
    $$
    \frac k2 + \frac{N(p-2)}{4} < 1.
    $$
    Up to a subsequence, we assume that $u_n \rightharpoonup u_c$ weakly in $H_0^1(B_1)$, $u_n \to u_c$ strongly in $L^2(B_1)$ and $L^p(B_1,|x|^{-k})$. It is clear that $u_c \in S_{2,c}$ and $E_2(u_c) \geq m_2(c)$. By the weak convergence of $u_n$ to $u_c$ in $H_0^1(B_1)$ and strong convergence in $L^p(B_1,|x|^{-k})$, we have
    \begin{align*}
    	E_2(u_c) \leq \lim_{n \to \infty}E_2(u_n) = m_2(c).
    \end{align*}
    Thus $E_2(u_c) = m_2(c)$ and $u_c$ is a normalized ground state of \eqref{eq1.2} on $S_{2,c}$.

	\vskip0.1in
	\emph{Step 2:} Any normalized ground state of \eqref{eq1.2} is positive after multiplying $-1$ if necessary and has Morse index $1$.

If $u_c$ satisfies $E_2(u_c) = \inf_{u \in S_{2,c}}E_2(u)$. Noticing that $E_2(|u_c|) \leq E_2(u_c)$, we know that $|u_c|$ is also a normalized ground state, which solves \eqref{eq1.2} for some $\lambda$. By the strong maximum principle, $|u_c| > 0$. Since $u_c$ is continuous, $u_c$ cannot change its sign. After multiplying $-1$ if necessary, we have $u_c = |u_c| > 0$.

Quiet similar to the proof of Lemma \ref{lem5.8}, we can prove that the Morse index of $u_c$ is $1$.

\vskip0.1in
\emph{Step 3:} If $u > 0$ solves \eqref{eq1.2} for some $\lambda$, then $\lambda < \lambda_1(B_1)$.

Let $e_1$ be the positive, unit eigenfunction corresponding to $\lambda_1(B_1)$. Then we have
\begin{equation}
\int_{B_1}|x|^{-k}|u|^{p-2}u_ce_1dx = \int_{B_1}(\nabla u \nabla e_1 - \lambda u e_1)dx = (\lambda_1(B_1) - \lambda)\int_{B_1} u e_1dx.
\end{equation}
Since $u > 0$, $\int_{B_1}|x|^{-k}|u|^{p-2}ue_1dx > 0$. Noticing that $\int_{B_1} u e_1dx > 0$, we know that $\lambda < \lambda_1$.

\vskip0.1in
\emph{Step 4:} Completion of the proof.

Suppose on the contrary that there are two positive normalized ground states $u_{c,1}, u_{c,2}$ with same $L^2$ mass $c$. Then $u_{c,1}, u_{c,2}$ solve \eqref{eq1.2} for some $\lambda, \widehat{\lambda}$ respectively. By Step 3, $\lambda < \lambda_1(B_1), \widehat{\lambda} < \lambda_1(B_1)$. By Step 2 and Theorem \ref{thm2.7}, both $u_{c,1}, u_{c,2}$ are on the global branch. Then, the fact $\int_{B_1}|u_{c,1}|^2dx = \int_{B_1}|u_{c,2}|^2dx$ implies that $\lambda = \widehat{\lambda}$ since the $L^2$ norm of solutions on the global branch is monotonic. Then we derive that $u_{c,1} = u_{c,2}$, completing the proof of Theorem \ref{thm2.8}.
\end{proof}

\subsection{Orbital stability results}

\begin{proof}[Proof of Theorem \ref{thm2.9}]
	Note that $u_\lambda$ is non-degenerate. Then expressing in our context the abstract theory developed in \cite{GSS1}, if $\partial_\lambda\int_{B_1}|u_\lambda|^2dx \neq 0$, then $\Psi(t,x) = e^{-i\lambda t}u_\lambda(x)$ is orbitally stable if $\partial_\lambda\int_{B_1}|u_\lambda|^2dx < 0$. Thus, thanks to Theorem \ref{thm2.7}, we can complete the proof.
\end{proof}

\appendix

\section{A counter-example where the normalized ground state exists but is not unique}\label{Aa}

Let $f(t) = t^{p}$ if $t \geq 0$ and $f(t) = t^{q}$ if $t < 0$ where $2 < p \neq q < 2 + 4s/N$. Let
\begin{align*}
	E(u) & = \frac{1}{2}\int_{\mathbb{R}^N}|(-\Delta)^{\frac{s}{2}} u|^2dx - \int_{\mathbb{R}^N}F(u)dx, \\
	Q(u) & = \frac{1}{2}\int_{\mathbb{R}^N}|u|^2dx, \\
	m(c) & = \inf_{u \in S_c}E(u) \quad \text{where } S_c = \big\{u \in H^s(\mathbb{R}^N): Q(u) = c\big\},
\end{align*}
where $F(t) = \int_0^tf(s)ds$. Consider the minimization problem
\begin{center}
	$(P_c)$      \ \ \ \ \ \ \  minimize $E(u)$ in $S_c$, i.e. find $u \in S_c$ such that $E(u) = m(c)$.
\end{center}

\begin{lemma}
	If $c = c_1 + c_2$ and $c_1, c_2 > 0$, then $m(c) < m(c_1) + m(c_2)$.
\end{lemma}

\begin{proof}
	The proof is similar to the one of Lemma \ref{lem 5.3}.
\end{proof}  

\begin{lemma}
	For any $c > 0$, there exists $u \in S_c$ such that $E(u) = m(c)$.
\end{lemma}

\begin{proof}
	The proof is similar to the one of Lemma \ref{lem5.2}.
\end{proof}

\begin{lemma}
	If $u \in S_c$ satisfies $E(u) = m(c)$, then $u \geq 0$ or $u \leq 0$.
\end{lemma}

\begin{proof}
	Let $c_\pm = Q(u_\pm)$ where $u_+ = \max\{u,0\}$, $u_- = \min\{u,0\}$. Arguing by contradiction, we assume that $c_\pm > 0$. Then $c = c_+ + c_-$ and $m(c) = m(c_+) + m(c_-) > m(c)$, which is a contradiction.
\end{proof}

\begin{corollary}
	Define
	\begin{align*}
		E^+(u) & = \frac{1}{2}\int_{\mathbb{R}^N}|(-\Delta)^{\frac{s}{2}} u|^2dx - \int_{\mathbb{R}^N}|u|^pdx, \quad m^+(c) = \inf_{u \in S_c}E^+(u), \\
		E^-(u) & = \frac{1}{2}\int_{\mathbb{R}^N}|(-\Delta)^{\frac{s}{2}} u|^2dx - \int_{\mathbb{R}^N}|u|^qdx, \quad m^-(c) = \inf_{u \in S_c}E^-(u).
	\end{align*}
	Then $m(c) =  \min\big\{m^+(c),m^-(c)\big\}$.
\end{corollary}

By \cite[Theorem 1.2 $(i)$]{Luo-Zhang}, $m^+(c) = c^{\frac{2ps}{4s-(p-2)N}}m^+(1)$, $m^-(c) = c^{\frac{2qs}{4s-(q-2)N}}m^-(1)$. Hence, there exists a unique $\hat{c} > 0$ such that $m^+(\hat{c}) = m^-(\hat{c})$. Then we know that the minimizers of $(P_c)$ are not unique at $c = \hat{c}$. Further, let $u_1$ and $u_2$ be minimizers for $m^+(\hat{c})$ and $m^-(\hat{c})$ respectively. Since $p \neq q$, it is clear that $u_1 \neq -u_2$.

\section{Existence and uniqueness of the global branch} \label{exi and uni}

Let $\Phi_{\lambda} = S(u) + G(u) - F(u) - \lambda Q(u)$. We will give results on existence and uniqueness of the global branch under a very general setting in this appendix. This has been done in \cite[Section 2]{HS} and we omit the details for the proofs here. We assume

\begin{itemize}
	\item[$(H_1)$] $W$ is a Hilbert space, $D_{u}Q(u) = u$, the linear operator $D_{u}S + D_{u}G$ is self-adjoint and bounded below on a Hilbert space $E$ with operator domain $\widehat{W}$ and form domain $W$.
\end{itemize}
Then $S(u) = \frac{1}{2}\left\langle D_{u}S(u),u\right\rangle, G(u) = \frac{1}{2}\left\langle D_{u}G(u),u\right\rangle, Q(u) = \frac{1}{2}\left\langle u,u\right\rangle$. Moreover, if $\lambda < \lambda_1 := \inf \sigma(D_{u}S + D_{u}G)$ where $\sigma(D_{u}S + D_{u}G)$ denotes the spectrum of $D_{u}S + D_{u}G: \widehat{W} \rightarrow E$, then
$$
\|u\|_{\lambda} = \sqrt{\left\langle D_{u}S(u) + D_{u}G(u) - \lambda u,u\right\rangle}
$$
is equivalent to the norm of $W$.

Let $G$ be a topological group, and assume that the action of $G$ on $E$ is isometric. Define
\begin{equation}
E_G = \big\{u \in E: gu = u, \forall g \in G\big\}, \quad \widehat{W}_G = \widehat{W} \cap E_G, \quad W_G = W \cap E_G.
\end{equation}
We assume the following:
\begin{itemize}
	\item[$(H_2)$] $\Phi_{\lambda}(gu) = \Phi_{\lambda}(u), \forall g \in G$.
	\item[$(H_3)$] $\Phi_{\lambda}|_{W_G}$ satisfies (PS) condition in $W_G$ for all $\lambda < \lambda_1$.
	\item[$(H_4)$] $D_{uu}F(0) = 0$, $\left\langle D_{uu}F(u)v,w \right\rangle = \left\langle v,D_{uu}F(u)w \right\rangle$. And there exists some $p > 2$ such that
	\begin{equation}
		\left\langle D_{uu}F(u)u,u \right\rangle \geq (p-1)\left\langle D_uF(u),u \right\rangle > 0, \forall u \in W \setminus \{0\}.
	\end{equation}
\end{itemize}

By the principle of symmetric criticality (see \cite[Theorem 1.28]{Wi}), $(H_2)$ ensures that any critical point of $\Phi_\lambda$ restricted to $\widehat{W}_G$ is a critical point of $\Phi_\lambda$.  To show the existence of solutions, we define
$$
\mathcal{N}_{G,\lambda} = \big\{u \in W_G \setminus \{0\}: \left\langle D_uS(u) + D_uG(u) - \lambda u,u \right\rangle = \left\langle D_uF(u) ,u \right\rangle\big\}.
$$
\begin{definition}
	Set $h_G(\lambda) = \inf_{u \in \mathcal{N}_{G,\lambda}}\Phi_\lambda(u)$. We say that $u \in W_G \setminus \{0\}$ is a $G$-ground state solution if $u$ solves (\ref{eqabs}) and achieves $h_G(\lambda)$.
\end{definition}

If $\mathcal{N}_{G,\lambda}$ is a $C^1$ manifold with codimension $1$ in $W_G$, $h_G(\lambda)$ is well-defined and $h_G(\lambda) \neq 0$ for all $\lambda < \lambda_1$, by $(H_3)$, standard arguments yield the existence of a $G$-ground state solution with $G$-Morse index $1$ (the definition will be given by \eqref{eqb.3}). Then, we establish a $C^1$ global branch in $\widehat{W}_G$ under the following assumptions:
\begin{itemize}
	\item[$(H_5)$] If $\Phi_\lambda'(u_\lambda) = 0$ and $\mu_{G}(u_\lambda) = 1$, then $\ker D_{uu}\Phi_\lambda(u_\lambda)|_{E_G} = \{0\}$, where $\mu_{G}(u)$, the $G$-Morse index of $u$, is defined as
	\begin{equation} \label{eqb.3}
		\mu_{G}(u):= \sharp \big\{e < 0: e \text{ is an eigenvalue of } D_{uu}\Phi_\lambda(u)|_{E_G}\big\}.
	\end{equation}
    \item[$(H_6)$] If $u \in W$ solves (\ref{eqabs}). Then $u \in \widehat{W}$.
    \item[$(H_7)$] $D_{u}F(u) \in E, \forall u \in W$, and $D_{uu}F(u)$ maps $\widehat{W}$ to $E$ for all $u \in \widehat{W}$.
\end{itemize}

Our results read as follows:

\begin{theorem}  \label{exi}
	Assume that $(H_1)-(H_7)$ hold, and that $\mathcal{N}_{G,\lambda}$ is a $C^1$ manifold with codimension $1$ in $W_G$, $h_G(\lambda)$ is well-defined and $h_G(\lambda) \neq 0$ for all $\lambda < \lambda_1$. Then there exists
	$$
	\lambda \mapsto u_\lambda \in C^1((-\infty,\lambda_1),W_G \setminus \{0\}),
	$$
	such that $u_\lambda$ solves (\ref{eqabs}), and $\mu_{G}(u_\lambda) = 1$.
\end{theorem}

\begin{theorem}[Uniqueness] \label{uni}
	Under the hypotheses of Theorem \ref{exi}. If we also assume that there exists $\Lambda_1 < \lambda_1$ such that (\ref{eqabs}) admits a unique solution with $G$-Morse index $1$ for all $\lambda < \Lambda_1$. Then for all $\lambda < \lambda_1$, (\ref{eqabs}) admits a unique solution with $G$-Morse index $1$.
\end{theorem}

\section{The relationship between the existence/uniqueness of the normalized ground state and the existence/uniqueness of ground state solution under Nehari manifold} \label{Neh}

In this appendix, we aim to show that the existence/uniqueness of normalized ground state under the constraint of $S_c$ imply these properties for the ground state solution under Nehari manifold. Define
\begin{align*}
	h(\lambda) := \inf_{u \in \mathcal{N}_{\lambda}}\Phi_\lambda(u) \quad \text{where } \mathcal{N}_{\lambda} = \big\{u \in W \setminus \{0\}: D_u\Phi_\lambda(u)(u) = 0\big\}.
\end{align*}
We say that $u \in W \setminus \{0\}$ is a ground state solution if $u$ achieves $h(\lambda)$, i.e. $\Phi_\lambda(u) = h(\lambda)$. Assume that $h(\lambda)$ is well-defined and $h(\lambda) \neq 0$ for all $\lambda < \lambda_1$. We have the following assumption:
\begin{itemize}
	\item[$(N_1)$] For any $u \in W \setminus \{0\}$, there exists a unique $t(u) = t(\lambda,u) > 0$ such that $t(u)u \in \mathcal{N}_{\lambda}$, and $\Phi_\lambda(t(u)u) = \max_{t > 0}\Phi_\lambda(tu)$.
\end{itemize}

\begin{lemma} \label{App c-lem1}
	Assume that $(N_1)$ holds. For any $c > 0$, $u \in \mathcal{N}_{\lambda}$, we have
	\begin{equation} \label{App c-1}
	\Phi_\lambda(u) \geq m(c) - \lambda c.
	\end{equation}
	Furthermore, the "=" holds if and only if $u$ is a normalized ground state on $S_c$ ($Q(u) = c$ and $E(u) = m(c)$), and $u$ is a ground state on $\mathcal{N}_{\lambda}$.
\end{lemma}

\begin{proof}
	Let $k = Q(u)$. Since $u \in \mathcal{N}_{\lambda}$, by $(N_1)$, one gets that $\Phi_\lambda(u) = \max_{t > 0}\Phi_\lambda(tu) \geq \Phi_\lambda(tu)$, and $\Phi_\lambda(u) = \Phi_\lambda(tu)$ if and only if $t = 1$. Then, by the definition of $m(c)$, we have
	\begin{eqnarray}
		\Phi_\lambda(u) \geq \Phi_\lambda(\sqrt{\frac{c}{k}}u) = E(\sqrt{\frac{c}{k}}u) - \lambda c \geq m(c) - \lambda c.
	\end{eqnarray}
	Thus \eqref{App c-1} holds true. On the one hand, if the "=" holds, then $E(\sqrt{\frac{c}{k}}u) = m(c)$ and $\Phi_\lambda(u) = \Phi_\lambda(\sqrt{\frac{c}{k}}u)$. By $(N_1)$, the latter implies that $k = c$, i.e. $Q(u)=c$. Hence, $u$ is a normalized ground state on $S_c$. Then by \eqref{App c-1}, for any $v \in \mathcal{N}_{\lambda}$, we have
	\begin{equation}
		\Phi_\lambda(v) \geq m(c) - \lambda c = E(u) - \lambda Q(u) = \Phi_\lambda(u).
	\end{equation}
	Thus $u$ is a ground state on $\mathcal{N}_{\lambda}$. On the other hand, if $u$ is a normalized ground state on $S_c$, one can deduce that $\Phi_\lambda(u) = E(u) - \lambda Q(u) = m(c) - \lambda c$. We know that the "=" holds. The proof is complete.
\end{proof}

\begin{remark}
	\eqref{App c-1} was used in \cite{DST} to study the equation
	\begin{equation}
	\Delta u + |u|^{p-2}u + \lambda u = 0 \quad \text{in } \ \Omega \subset \mathbb{R}^N.
	\end{equation}
\end{remark}

\begin{theorem} \label{App c-thm3}
	Assume that $(N_1)$ holds. If $u$ is a normalized ground state on $S_c$, then $\Phi_\lambda(u) = h(\lambda)$, where $\lambda$ is the Lagrange multiplier of $u$. In particular, the existence of normalized ground state under the constraint of $S_c$ implies the one of ground state solution under Nehari manifold.
\end{theorem}

\begin{proof}
	Let $c = Q(u)$ and $\lambda$ be the Lagrange multiplier of $u$. By \eqref{App c-1}, for any $v \in \mathcal{N}_{\lambda}$, we have
	\begin{equation}
		\Phi_\lambda(v) \geq m(c) - \lambda c = E(u) - \lambda Q(u) = \Phi_\lambda(u),
	\end{equation}
	implying that $\Phi_\lambda(u) = h(\lambda)$.
\end{proof}

\begin{theorem}
	Assume that $(N_1)$ holds, the normalized ground state on $S_c$ is unique and is $u$. Let $\lambda$ be the Lagrange multiplier of $u$. Then the ground state solution on $\mathcal{N}_{\lambda}$ is unique.
\end{theorem}

\begin{proof}
	By Theorem \ref{App c-thm3}, $u$ is a ground state solution on $\mathcal{N}_{\lambda}$.  Let us argue by contradiction. Assume that there exists another ground state $v \in \mathcal{N}_{\lambda}$, which satisfies $\Phi_\lambda(v) = h(\lambda)$. Take $c = Q(u)$ and we have
	\begin{equation}
		\Phi_\lambda(v) = \Phi_\lambda(u) = m(c) - \lambda c.
	\end{equation}
	For $v \in \mathcal{N}_{\lambda}$, the "=" in \eqref{App c-1} holds. Then by Lemma \ref{App c-lem1},  $v$ is a normalized ground state on $S_c$, which contradicts with the uniqueness of $u$. The proof is complete.
\end{proof}

\begin{remark}
	For all $c > 0$, assume that the normalized ground state $u_c$ on $S_c$ is unique and $\lambda(c)$ is the Lagrange multiplier of $u_c$. In concrete applications, it can be proved that $\lambda(c)$ is non-increasing with respect to $c > 0$ (see the proof of Theorem \ref{thm2.10}). By the uniqueness of $u_c$ for all $c > 0$, it is easy to verify that $\lambda(c)$ is strictly decreasing on $c > 0$. Then, also by the uniqueness of $u_c$, $\lambda(c)$ is continuous on $c > 0$. If $(\lim_{c \rightarrow \infty}\lambda(c),\lim_{c \rightarrow 0^+}\lambda(c)) = (-\infty,\lambda_1)$, then the ground state solution on $\mathcal{N}_{\lambda}$ is unique for all $\lambda < \lambda_1$.
\end{remark}

\section{The expression of $m(c)$}

Take \eqref{eq1.4} as an example in this appendix. If the normalized ground state $u_c$ is unique for any $c > 0$, Theorem \ref{thm2.10} $(ii)$ shows us that $m'(c) = \lambda(c)$ where $\lambda(c)$ is the Lagrange multiplier of $u_c$. By the Pohozaev identity, we have
$$
m'(c) = \lambda(c) = \frac{1}{c}\int_{\mathbb{R}^N}\left( \frac{N-2s}{N}|(-\Delta)^{\frac{s}{2}} u_c|^2 - \frac{2}{p}h(|x|)|u_c|^{p} - \frac{2}{pN}h'(|x|)|x||u_c|^p\right)dx.
$$
Then $m(c)$ is the integral of $m'(c)$.

Here we provide another method to show the expression of $m(c)$. Let us set $s=1$ for the simplicity. We define $v(x)=tu(\beta x)$.   Then
\begin{eqnarray*}
	\|\nabla v\|_{L^2}^2=t^2\beta^{2-N}\|\nabla u\|_{L^2}^2\quad\text{and}\quad\|v\|_{L^p}^p=t^p\beta^{-N}\|u\|_{L^p}^p.
\end{eqnarray*}
Now, by solving the system
\begin{eqnarray*}
	\left\{\aligned&t^2\beta^{2-N}=t^p\beta^{-N},\\
	&t^2\beta^{-N}=(1\pm\frac{\varepsilon}{c})^2,\endaligned\right.
\end{eqnarray*}
we have
\begin{eqnarray*}
	t_{\varepsilon,\pm}=\beta_{\varepsilon,\pm}^{\frac{2}{p-2}}\quad\text{and}\quad \beta_{\varepsilon,\pm}=(1\pm\frac{\varepsilon}{c})^{\frac{2(p-2)}{4+2N-Np}},
\end{eqnarray*}
which implies that
\begin{eqnarray*}
	\beta_{\varepsilon,\pm}=1\pm\frac{2(p-2)\varepsilon}{(4+2N-Np)c}+O(\varepsilon^2).
\end{eqnarray*}
Let $u_c$ be the minimizer of $m(c)$.  Then $v_{\varepsilon,\pm}\in S_{c\pm\varepsilon}$, where $v_{\varepsilon,\pm}(x)=t_{\varepsilon,\pm} u_c(\beta_{\varepsilon,\pm} x)$.
Thus,
\begin{eqnarray*}
	&& m(c\pm\varepsilon)\\
	&\leq&\frac{1}{2}\|\nabla v_{\varepsilon,\pm}\|_{L^2}^2-\frac{1}{p}\int_{\mathbb{R}^N}h(|x|)|v_{\varepsilon,\pm}|^pdx\\
	&=&\beta_{\varepsilon,\pm}^{\frac{4}{p-2}-N+2}(m(c)-\frac1p\int_{\mathbb{R}^N}(h(\beta_{\varepsilon,\pm}^{-1}|x|)-h(|x|))|u_c|^pdx)\\
	&=&(1\pm\frac{\varepsilon}{c})^{\frac{2(2N-(N-2)p)}{4+2N-Np}}(m(c)-\frac1p\int_{\mathbb{R}^N}(h(\beta_{\varepsilon,\pm}^{-1}|x|)-h(|x|))|u_c|^pdx)\\
	&=&(1\pm\frac{2(2N-(N-2)p)\varepsilon}{(4+2N-Np)c}+O(\varepsilon^2))(m(c)\\
	&& \pm\frac{2(p-2)\varepsilon}{(4+2N-Np)pc}\int_{\mathbb{R}^N}h'(|x|)|x||u_c|^pdx+O(\varepsilon^2))\\
	&=&m(c)\pm(\frac{2(2N-(N-2)p)m(c)}{(4+2N-Np)c}+\frac{2(p-2)}{(4+2N-Np)pc}\int_{\mathbb{R}^N}h'(|x|)|x||u_c|^pdx)\varepsilon+O(\varepsilon^2)
\end{eqnarray*}
if $h(|x|)$ has the good regularity, such as $h(r)\in W^{2,\frac{2^*}{2^*-p}}([0, +\infty))$ since $u_c$ decays to zero as $|x|\to+\infty$, exponentially.  It follows that
\begin{eqnarray*}
	m'_{+}(c)&=&\lim_{\varepsilon\to0}\frac{m(c+\varepsilon)-m(c)}{\varepsilon}\\
	&\leq&\frac{2(2N-(N-2)p)m(c)}{(4+2N-Np)c}+\frac{2(p-2)}{(4+2N-Np)pc}\int_{\mathbb{R}^N}h'(|x|)|x||u_c|^pdx
\end{eqnarray*}
and
\begin{eqnarray*}
	m'_{-}(c)&=&\lim_{\varepsilon\to0}\frac{m(c-\varepsilon)-m(c)}{-\varepsilon}\\
	&\geq&\frac{2(2N-(N-2)p)m(c)}{(4+2N-Np)c}+\frac{2(p-2)}{(4+2N-Np)pc}\int_{\mathbb{R}^N}h'(|x|)|x||u_c|^pdx
\end{eqnarray*}
which implies that $m'_{+}(c)\leq m'_{-}(c)$.

We do not need to assume that $u_c$ is unique in the above discussions. However, from now on, let $u_c$ be unique. On the other hand, by solving the system
\begin{eqnarray*}
	\left\{\aligned&t^2\beta^{2-N}=t^p\beta^{-N},\\
	&t^2\beta^{-N}=(\frac{c}{c\pm\varepsilon})^2,\endaligned\right.
\end{eqnarray*}
we have
\begin{eqnarray*}
	\widetilde{t}_{\varepsilon,\pm}=\widetilde{\beta}_{\varepsilon,\pm}^{\frac{2}{p-2}}\quad\text{and}\quad \widetilde{\beta}_{\varepsilon,\pm}=(\frac{c}{c\pm\varepsilon})^{\frac{2(p-2)}{4+2N-Np}},
\end{eqnarray*}
	which implies that
\begin{eqnarray*}
	\widetilde{\beta}_{\varepsilon,\pm}=1\mp\frac{2(p-2)\varepsilon}{(4+2N-Np)c}+O(\varepsilon^2).
\end{eqnarray*}
Now, let $u_{c\pm\varepsilon}$ be the minimizer of $m(c\pm\varepsilon)$.  Then $\widetilde{v}_{\varepsilon,\pm}\in S_{c}$, where $\widetilde{v}_{\varepsilon,\pm}(x)=\widetilde{t}_{\varepsilon,\pm} u_{c\pm\varepsilon}(\widetilde{\beta}_{\varepsilon,\pm} x)$.  Thus,
\begin{eqnarray*}
	&& m(c)\\
	&\leq&\frac{1}{2}\|\nabla \widetilde{v}_{\varepsilon,\pm}\|_{L^2}^2-\frac{1}{p}\int_{\mathbb{R}^N}h(|x|)|\widetilde{v}_{\varepsilon,\pm}|^pdx\\
	&=&\widetilde{\beta}_{\varepsilon,\pm}^{\frac{4}{p-2}-N+2}(m(c\pm\varepsilon)-\frac1p\int_{\mathbb{R}^N}(h(\widetilde{\beta}_{\varepsilon,\pm}^{-1}|x|)-h(|x|))|u_{c\pm\varepsilon}|^pdx)\\
	&=&(1\mp\frac{\varepsilon}{c}+O(\varepsilon^2))^{\frac{2(2N-(N-2)p)}{4+2N-Np}}(m(c\pm\varepsilon)
	-\frac1p\int_{\mathbb{R}^N}(h(\widetilde{\beta}_{\varepsilon,\pm}^{-1}|x|)-h(|x|))|u_{c\pm\varepsilon}|^pdx)\\
	&=&(1\mp\frac{2(2N-(N-2)p)\varepsilon}{(4+2N-Np)c}+O(\varepsilon^2))(m(c\pm\varepsilon)\\
	&& \mp\frac{2(p-2)\varepsilon}{(4+2N-Np)pc}\int_{\mathbb{R}^N}h'(|x|)|x||u_{c\pm\varepsilon}|^pdx+O(\varepsilon^2))\\
	&=& m(c\pm\varepsilon)\mp(\frac{2(2N-(N-2)p)m(c\pm\varepsilon)}{(4+2N-Np)c}\\
	&& +\frac{2(p-2)}{(4+2N-Np)pc}\int_{\mathbb{R}^N}h'(|x|)|x||u_{c\pm\varepsilon}|^pdx)\varepsilon+O(\varepsilon^2).
\end{eqnarray*}
It follows from the continuity of $m(c)$ that
\begin{eqnarray*}
	m'_{+}(c)&=&\lim_{\varepsilon\to0}\frac{m(c+\varepsilon)-m(c)}{\varepsilon}\\
	&\geq&\frac{2(2N-(N-2)p)m(c)}{(4+2N-Np)c}+\frac{2(p-2)}{(4+2N-Np)pc}\int_{\mathbb{R}^N}h'(|x|)|x||u_c|^pdx
\end{eqnarray*}
and
\begin{eqnarray*}
	m'_{-}(c)&=&\lim_{\varepsilon\to0}\frac{m(c-\varepsilon)-m(c)}{-\varepsilon}\\
	&\leq&\frac{2(2N-(N-2)p)m(c)}{(4+2N-Np)c}+\frac{2(p-2)}{(4+2N-Np)pc}\int_{\mathbb{R}^N}h'(|x|)|x||u_c|^pdx
\end{eqnarray*}
which implies that $m'_{+}(c)\geq m'_{-}(c)$.  Thus, we always have $m'_{+}(c)= m'_{-}(c)$ which implies that $m'(c)$ exists for all $c>0$.  Moreover, $u_c$ exists only for $2<p<2+\frac{4}{N}$.  Thus, by
\begin{eqnarray*}
	m'(c)=\frac{2(2N-(N-2)p)m(c)}{(4+2N-Np)c}+\frac{2(p-2)}{(4+2N-Np)pc}\int_{\mathbb{R}^N}h'(|x|)|x||u_c|^pdx,
\end{eqnarray*}
we have
\begin{eqnarray*}
	m(c)=c^{-\frac{2(2N-(N-2)p)}{(4+2N-Np)}}\left(m(1)+\frac{2(p-2)}{(4+2N-Np)p}\int_{1}^{c}\left(\int_{\mathbb{R}^N}h'(|x|)|x||u_\tau|^pdx\right)d\tau\right).
\end{eqnarray*}
															
For general nonlinearities $F(|x|,u)$, the above calculations are still valid if $F(|x|,u)$ also have the good regularity, such as $\partial_{rr}F(r,u)r^2\in L^{1}(\mathbb{R}^N)$.\\

\

\textbf{Acknowledgements}

\
This work was started when Song was a Phd student at Institute of Mathematics, AMSS, Chinese Academy of Science, and revised when he was visiting the Universit\'e Marie et Louis Pasteur, LmB (UMR 6623). 
The authors thank the anonymous referees so much for careful reviewing this article and providing valuable revision suggestions.

\

\textbf{Conflict of Interest Statement}

\

On behalf of all authors, the corresponding author states that there is no conflict of interest.

\

\

\textbf{Data Availability Statement}

\

My manuscript has no associated data.

\end{document}